\documentclass[11pt,oneside,reqno]{amsart}

\usepackage{amssymb}
\usepackage{algpseudocode}
\usepackage{algorithm}
\usepackage{verbatim}
\usepackage{graphicx}
\usepackage[noadjust]{cite}
\RequirePackage{dsfont} 
 \scrollmode
\allowdisplaybreaks
\usepackage{graphicx}
\usepackage{epstopdf}
\usepackage[usenames]{color}
\definecolor{darkred}{RGB}{139,0,0}
\definecolor{darkgreen}{RGB}{0,100,0}
\definecolor{darkmagenta}{RGB}{139,0,139}

\usepackage{amsmath}
\usepackage{tikz}

\makeatletter
\newcommand{\xleftrightarrow}[2][]{\ext@arrow 3359\leftrightarrowfill@{#1}{#2}}
\makeatother

\newcommand{\xdasharrow}[2][->]{
\tikz[baseline=-\the\dimexpr\fontdimen22\textfont2\relax]{
\node[anchor=south,font=\scriptsize, inner ysep=1.5pt,outer xsep=2.2pt](x){#2};
\draw[shorten <=3.4pt,shorten >=3.4pt,dashed,#1](x.south west)--(x.south east);
}
}

\newcommand{\DEBUG}{}

\ifdefined\DEBUG%

  \def\rem#1{{\marginpar{\raggedright\scriptsize #1}}}
  \newcommand{\pmr}[1]{\rem{\color{blue}{$\bullet$ #1}}}
  \newcommand{\ppr}[1]{\rem{\color{red}{$\bullet$ #1}}}
 \else%

  \newcommand{\ppr}[1]{}
  \newcommand{\pmr}[1]{}
 \fi

\def\rho{\varrho}

\def\rd{\,{\mathrm d}}

\theoremstyle{plain}
\newtheorem{theorem}{Theorem}
\newtheorem{lemma}{Lemma}
\newtheorem{fact}{Fact}
\newtheorem{corollary}{Corollary}
\newtheorem{proposition}{Proposition}

\theoremstyle{definition}
\newtheorem{remark}{Remark}

\begin{document}

\title
[Randomized Runge-Kutta method]
{Randomized Runge-Kutta method -- stability and convergence under inexact information}

\author[T. Bochacik]{Tomasz Bochacik}
\address{AGH University of Science and Technology,
Faculty of Applied Mathematics,
Al. A.~Mickiewicza 30, 30-059 Krak\'ow, Poland}
\email{bochacik@agh.edu.pl}

\author[M. Goćwin]{Maciej Go\'cwin}
\address{AGH University of Science and Technology,
Faculty of Applied Mathematics,
Al. A.~Mickiewicza 30, 30-059 Krak\'ow, Poland}
\email{gocwin@agh.edu.pl}

\author[P. Morkisz]{Pawe{\l} M. Morkisz}
\address{AGH University of Science and Technology,
Faculty of Applied Mathematics,
Al. A.~Mickiewicza 30, 30-059 Krak\'ow, Poland}
\email{morkiszp@agh.edu.pl}

\author[P. Przyby{\l}owicz]{Pawe{\l} Przyby{\l}owicz}
\address{AGH University of Science and Technology,
Faculty of Applied Mathematics,
Al. A.~Mickiewicza 30, 30-059 Krak\'ow, Poland}
\email{pprzybyl@agh.edu.pl, corresponding author}

\begin{abstract}
We deal with optimal approximation of solutions of ODEs  under  local Lipschitz condition and inexact discrete information about the right-hand side functions. We show that the randomized two-stage  Runge-Kutta scheme is the optimal method among all randomized algorithms based on standard noisy information. We perform numerical experiments that confirm our theoretical findings. Moreover, for the optimal algorithm we rigorously investigate properties of regions of absolute stability. 
\end{abstract}
\maketitle
\tableofcontents
This paper is devoted to  the problem of optimal approximation of solutions of ordinary differential equations (ODEs) of the following form
\begin{equation}
	\label{PROBLEM}
		\left\{ \begin{array}{ll}
			z'(t)= f(t,z(t)), \ t\in [a,b], \\
			z(a) = \eta, 
		\end{array}\right.
\end{equation}
where $-\infty<a<b<\infty$, $d\in\mathbb{N}$, $\eta\in\mathbb{R}^d$, $f: [a,b]\times \mathbb{R}^d\to \mathbb{R}^d$.  We consider the case when $f$ is only locally Lipschitz. Due to the low regularity of the problem, we focus on the class of randomized algorithms. Moreover, we assume that we have access to $f$ only through its noisy evaluations. We aim at defining an algorithm that approximates $z$ optimally, i.e. with the minimal possible error. Moreover, we want to investigate stability properties of the optimal scheme.

Approximation of solutions of ODEs via randomized algorithms and under exact information about right-hand side functions is a problem well studied in the literature, see, for example,\cite{daun1, hein1, HeinMilla, JenNeuen, Kac1, KruseWu_1, stengle1, stengle2}. However, there are still few papers on  approximate solving (even via deterministic algorithms) of ODEs when the available information is corrupted, see \cite{KaPr16}. Inexact information has been mainly investigated in the context of function integration and approximation (\cite{MoPl16, MoPl20}), approximate solving of PDEs (\cite{Wer96, Wer97}), stochastic integration and SDEs (\cite{KaMoPr19, PMPP17, PMPP19}). Such analysis, under noisy information, seems to be important from the  point of view of applications and stability issues, see \cite{Pla96} and Remarks \ref{GPU_model}, \ref{zero_stab_RK}.  We also refer to \cite{JoC4, KaPl90, milvic} for further discussion and examples.

In this paper we extend the results concerning randomized Runge-Kutta scheme (known from \cite{KruseWu_1}) in three directions. Firstly, we investigate the error and optimality of the randomized Runge-Kutta method in the case when $f$ is only locally Lipschitz. Secondly, we allow noisy evaluations of $f$. This means that the (randomized) evaluations of $f$ might be corrupted by some noise at level of at most $\delta\in [0,1]$, which corresponds to the precision level. Finally, we rigorously prove fundamental properties of regions of  stability, such as openness, boundedness, symmetry. We consider three types of such regions due to the three types of convergence of underlying sequences of random variables: mean-square, with probability $1$, and in probability. For the stability analysis we adopt the approach used in \cite{higham2000} in the context of stochastic differential equations.

The novelty and main results of the paper can be summarized as follows:
\begin{itemize}
    \item We present upper bound on the $L^p(\Omega)$-error for the randomized Runge-Kutta method in the presence of informational noise and under local Lipschitz condition (Theorem \ref{error_RRK}). We emphasise a strong connection between the error analysis  under inexact information and $0$-stability of the method (Remark \ref{zero_stab_RK}).
    \item  We show respective lower bound and then we justify that the randomized Runge-Kutta  scheme is optimal in the class of locally Lipschitz right-hand side functions $f$, among all algorithms based on randomized inexact  information  about $f$  (Theorem \ref{thm_opt_RK}).
    \item We rigorously prove  properties  of  regions of stability for the randomized Runge-Kutta scheme (Theorems \ref{prop_MS}, \ref{prop_AS} and the equality \eqref{RSP}).  According to our best knowledge this is a first attempt in this direction.
\end{itemize}
The paper is organized into seven sections. Section 1 contains problem definition and description of the used model of computation under noisy information. Upper bounds on the error of randomized Runge-Kutta methods are established in Section 2. Corresponding lower bound and optimality are discussed in Section 3. In Section 4 we report the results of numerical experiments performed for two exemplary equations, where one of them is the SIR model. Properties of  regions of stability for the randomized Runge-Kutta method are investigated in Section 6. Finally, Appendix consists of some auxiliary results that are used in the paper. 
\section{Preliminaries}
Let $\|\cdot\|$ be the first norm in $\mathbb{R}^d$, i.e. $\|x\|=\sum\limits_{k=1}^d|x_k|$ for $x\in \mathbb{R}^d$. By $\{e_k\}_{k=1}^d$ we denote the canonical base in $\mathbb{R}^d$. For $x\in\mathbb{R}^d$ and $r\in [0,\infty)$ we denote by $B(x,r)=\{y\in\mathbb{R}^d \ | \ \|y-x\|\leq r\}$ the closed ball in $\mathbb{R}^d$. Moreover, we write $\mathbb{C}_{-}=\{z\in\mathbb{C} \colon  \Re(z)<0\}$. Let $(\Omega,\Sigma,\mathbb{P})$ be a complete probability space. For a random variable $X:\Omega\to\mathbb{R}$, defined on  $(\Omega,\Sigma,\mathbb{P})$, we denote by $\|X\|_p=(\mathbb{E}|X|^p)^{1/p}$, $p\in [2,\infty)$. For a Polish space $E$ by $\mathcal{B}(E)$ we denote the Borel $\sigma$-field on $E$. 

Let $K,L\in(0,\infty)$, $\varrho\in(0,1]$. We will consider a class $F^{\varrho}=F^{\varrho}(K,L)$ of pairs $(\eta,f)$ defined by the following conditions:
\begin{itemize}
	\item [(A0)] $\|\eta\| \leq K$,
	\item [(A1)] $f\in C\left([a,b]\times\mathbb{R}^d\right)$,
	\item [(A2)] $\|f(t,x)\| \leq K (1+\|x\|),$ for all $(t,x)\in [a,b]\times\mathbb{R}^d$. 
\end{itemize}
Take $\bar R=\bar R(a,b,K)$ as 
\begin{eqnarray}
\label{def_bar_R}
    &&\bar R = \max\Bigl\{K(1+b-a)\Bigl(1+e^{K(b-a)}(1+K(b-a))\Bigr),\notag\\
    && K+(b-a)(1+K)+\Bigl(\frac{1}{K}+1\Bigr)(1+K(b-a))\Bigl(e^{K(b-a)(1+K(b-a))}(1+K)-1\Bigr)\Bigr\}.
\end{eqnarray}
By  Lemma \ref{Radii_R} below we will see that it is sufficient for our analysis to assume that $f$ satisfies Lipschitz condition only in the ball $B(\eta,\bar R)$.
Namely, in addition to $(A0),(A1),(A2)$ we assume that
\begin{itemize}
    \item [(A3)] $\|f(t,x)-f(s,x)\|\leq L|t-s|^{\varrho}$ for all $t,s\in [a,b]$, $x\in B(\eta,\bar R)$,
    \item [(A4)] $\|f(t,x)-f(t,y)\|\leq L\|x-y\|$ for all $t\in [a,b]$, $x,y\in B(\eta,\bar R)$.
\end{itemize}
The numbers $a,b,d,\varrho,K,L$ will be called parameters of the class $F^{\varrho}$. Except for $a$, $b$, and $d$ the parameters are, in general, not known and the algorithms presented later on will not use them as input parameters.

We wish to approximate solution of \eqref{PROBLEM}  for  $f\in F^{\varrho}$ by an algorithm that is based on noisy information about $f$. Namely, we assume that access to the function $f$ is possible only through its noisy evaluations
\begin{equation}
    \tilde f(t,y) = f(t,y) + \tilde \delta_f(t, y), \ (t,y)\in [a,b]\times\mathbb{R}^d,
\end{equation}
 where $\tilde \delta_f(t, y)$ is an error function corrupting the exact value $f(t,y)$, such that $ \|\tilde \delta_f(t, y) \|\leq \delta$.  We refer to $\delta$ as to precision parameter. Moreover, we allow randomized choice of the evaluation points $(t,y)$. We now describe model of computation in all details.
 
In order to define  model of computation under randomized inexact
information about $f$ we need to introduce the following auxiliary class
\begin{displaymath}
  \mathcal{K}(\delta)=\{\tilde \delta:[a,b]\times\mathbb{R}^d\to \mathbb{R}^d \colon \tilde\delta-\hbox{Borel measurable}, \
   \|\tilde\delta(t,y)\|\leq\delta \ \hbox{for all} \ t\in [a,b], y\in \mathbb{R}^d \},
\end{displaymath}
where we assume  for the  precision parameter $\delta$ that $\delta\in[0,1]$. Note that the constant mappings $\tilde\delta (t,y)=\pm \delta e_1$ belong to $\mathcal{K}(\delta)$. (This is important fact when establishing lower bounds, see \cite{PMPP17} and Section 4 below.) Moreover, for $(\eta,f)\in F^{\varrho}$ let
\begin{equation}
    V_{(\eta,f)}(\delta) = B(\eta,\delta)\times V_{f}(\delta),
\end{equation}
where
\begin{equation}
    V_{f}(\delta)=\{ \tilde f \colon \exists_{\tilde\delta_f\in\mathcal{K}(\delta)} \ \tilde f =f+\tilde\delta_f\}.
\end{equation}
It holds that $V_{(\eta,f)}(\delta)\subset V_{(\eta,f)}(\delta')$ for $0\leq\delta\leq\delta'\leq 1$ and $V_{(\eta,f)}(0)=\{\eta\}\times\{f\}$. 

For $(\eta,f)\in F^{\varrho}$ let $(\tilde\eta,\tilde f)\in V_{(\eta,f)}(\delta)$. A vector of noisy information about $(\eta,f)$ is as follows
\begin{equation}
    \mathcal{N}(\tilde \eta,\tilde f)=[\tilde f(t_0,y_0),\ldots,\tilde f(t_{i-1},y_{i-1}),\tilde f(\theta_0,z_0),\ldots,\tilde f(\theta_{i-1},z_{i-1}),\tilde\eta],
\end{equation}
where $i\in\mathbb{N}$ and $(\theta_0,\theta_1,\ldots,\theta_{i-1})$ is a random vector on  $(\Omega,\Sigma,\mathbb{P})$. For  Borel measurable mappings $\psi_j:\mathbb{R}^{(2j+1)d}\to\mathbb{R}^d\times \mathbb{R}^d$, $j\in\left\{0,\ldots,i-1\right\}$, we set
\begin{equation}
    (y_0,z_0)=\psi_0(\tilde\eta),
\end{equation}
and
\begin{equation}
    (y_j,z_j)=\psi_j\Bigl(\tilde f(t_0,y_0),\ldots,\tilde f(t_{j-1},y_{j-1}),\tilde f(\theta_0,z_0),\ldots,\tilde f(\theta_{j-1},z_{j-1}),\tilde\eta\Bigr)
\end{equation}
for $j\in\left\{1,\ldots,i-1\right\}$. The total number of noisy valuations of $f$ is $l=2i$. Note that $\mathcal{N}(\tilde \eta,\tilde f):\Omega\to\mathbb{R}^{(2i+1)d}$ is a random vector. 

Any algorithm $\mathcal{A}$ using $\mathcal{N}(\tilde\eta,\tilde f)$ that computes the approximation to $z$ is given by
\begin{equation}
\label{def_alg}
    \mathcal{A}(\tilde\eta,\tilde f,\delta)=\varphi(\mathcal{N}(\tilde \eta,\tilde f)),
\end{equation}
where
\begin{equation}
    \varphi:\mathbb{R}^{(2i+1)d}\to D([a,b];\mathbb{R}^d)
\end{equation}
is a Borel measurable function. In the Skorokhod space $D([a,b];\mathbb{R}^d)$ we consider the Borel $\sigma$-field $\mathcal{B}(D([a,b];\mathbb{R}^d))$ that coincides with the $\sigma$-field generated by coordinate mappings, see Theorem 7.1 in \cite{Parth}. This assures that $\mathcal{A}(\tilde\eta,\tilde f,\delta):\Omega\to D([a,b];\mathbb{R}^d)$ is $\Sigma$-to-$\mathcal{B}(D([a,b];\mathbb{R}^d))$ measurable and, by  Theorem  7.1 in \cite{Parth}, for all $t\in [a,b]$ the mapping
\begin{equation}
\label{fix_t_alg}
    \Omega\ni \omega\mapsto \mathcal{A}(\tilde\eta,\tilde f,\delta)(\omega)(t)\in\mathbb{R}^d
\end{equation}
is $\Sigma$-to-$\mathcal{B}(\mathbb{R}^d)$-measurable. For a given $n\in\mathbb{N}$ we denote by $\Phi_n$ a class of all algorithms of the form \eqref{def_alg} for which the total number of evaluations $l$ is at most $n$.

Let $p\in [2,\infty)$. For a fixed $(\eta,f)\in F^{\varrho}$ the error of $\mathcal{A}\in\Phi_n$ is given as
\begin{equation}
\label{fixed_err}
    e^{(p)}(\mathcal{A},\eta,f,\delta)=\sup\limits_{(\tilde \eta,\tilde f)\in V_{(\eta,f)}(\delta)} \Bigl\|\sup\limits_{a\leq t\leq b}\|z(\eta,f)(t)-\mathcal{A}(\tilde\eta,\tilde f,\delta)(t)\|\Bigl\|_p,
\end{equation}
see Remark \ref{well_def_err}. The worst-case error of the algorithm $\mathcal{A}$ is
defined by
\begin{equation}
    e^{(p)}(\mathcal{A},\mathcal{G},\delta)=\sup\limits_{(\eta,f)\in \mathcal{G}} e^{(p)}(\mathcal{A},\eta,f,\delta),
\end{equation}
where $\mathcal{G}$ is a certain subclass of $F^{\varrho}$, see \cite{TWW88}. Finally,
we consider the $n$th minimal error defined as follows
\begin{equation}
\label{nth_min_err}
    e^{(p)}_n(\mathcal{G},\delta)=\inf\limits_{\mathcal{A}\in\Phi_n}e^{(p)}(\mathcal{A},\mathcal{G},\delta).
\end{equation}
Our aim is two fold:
\begin{itemize}
    \item determine sharp bounds on the $n$th minimal error and to define implementable  algorithm for which the infimum in \eqref{nth_min_err} is asymptotically attained. We call such an  algorithm the optimal one, see \cite{TWW88}.
    \item investigate  stability (in certain sense) of the defined  optimal method.
\end{itemize}
We follow the usual convention that all constants appearing in this paper (including those in the "O", "$\Omega$", and "$\Theta$" notation) will only depend on the parameters of the  class $F^{\varrho}$ and $p$. Furthermore, the same symbol may be used for different constants.
\begin{remark}
\label{GPU_model}  We want to underline here that  proposed model of computation covers the phenomenon of lowering precision of computations. For example, in the scalar case (i.e. $d=1$)  we can model relative round off errors by considering the following disturbing functions $\tilde \delta_f$:
\begin{equation}
	\tilde\delta_f(t,y)=\delta\cdot\alpha(t,y)\cdot f(t,y), \ (t,y)\in [0,T]\times\mathbb{R},
\end{equation}
where $\alpha$ is a Borel measurable bounded function on $[0,T]\times\mathbb{R}$. This is a frequent case for efficient computations using both CPUs and GPUs. See \cite{KaMoPr19}, where similar model of noisy information was considered and Monte Carlo simulations were performed on GPUs. Moreover, in \cite{PMPP17, PMPP19} the authors show results of numerical experiments (performed on CPUs) concerning approximate solving of SDEs under inexact information.
\end{remark}
\begin{remark}
\label{well_def_err}
    Due to Lemma \ref{Radii_R} (i)  for all $\omega\in\Omega$ the mapping
$[a,b]\ni t\mapsto \|z(\eta,f)(t)-\mathcal{A}(\tilde\eta,\tilde f,\delta)(\omega)(t)\|\in [0,\infty)$ belongs to $D([a,b];[0,\infty))$, and by \eqref{fix_t_alg} for all $t\in [a,b]$ the function $\Omega\ni\omega\mapsto\|z(\eta,f)(t)-\mathcal{A}(\tilde\eta,\tilde f,\delta)(\omega)(t)\|\in [0,\infty)$
is a $\Sigma$-to-$\mathcal{B}([0,\infty))$ measurable. Hence,  $\Omega\ni\omega\mapsto\sup\limits_{a\leq t\leq b}\|z(\eta,f)(t)-\mathcal{A}(\tilde\eta,\tilde f,\delta)(\omega)(t)\|\in [0,\infty)$
is $\Sigma$-to-$\mathcal{B}([0,\infty))$ measurable function and the error \eqref{fixed_err} is well-defined.
\end{remark}
\section{Randomized Runge-Kutta method under noisy information}
In the case of inexact information about $(\eta,f)$  the randomized two-stage Runge-Kutta algorithm is defined as follows. Let $n\in\mathbb{N}$, $h=(b-a)/n$, and $t_j = a + jh$, $j\in\left\{0,\ldots,n\right\}$.
We assume that $(\tau_j)_{j\in\mathbb{N}}$ are independent and identically distributed random variables on  $(\Omega, \Sigma, \mathbb{P})$, uniformly distributed on $[0,1]$. Let $(\eta,f)\in F^{\varrho}$ and $(\tilde \eta, \tilde f)\in V_{(\eta,f)}(\delta)$. Then we set
\begin{equation}
\label{def_rand_RK_1}
		\left\{ \begin{array}{ll}
			\bar V^0 := \tilde\eta,\\
			\bar V_\tau^j := \bar V^{j-1} + h \tau_j \tilde f (t_{j-1},  \bar V^{j-1}), \\
			\bar V^j := \bar V^{j-1} + h \tilde f(\theta_j, \bar V_\tau^j)
		\end{array}\right.
\end{equation}
for $j\in\left\{1,\ldots,n\right\}$, where $\theta_j := t_{j-1} + \tau_j h$.
The final approximation of $z$ in the interval $[a,b]$ is obtained by taking
\begin{equation}
    \bar l(t)=\bar l_j(t), \ t\in [t_j,t_{j+1}], \ j\in\left\{0,\ldots,n-1\right\},
\end{equation}
where
\begin{equation}
\label{def_blj}
    \bar l_j(t) = \frac{\bar V^{j+1}-\bar V^{j}}{t_{j+1}-t_j}(t-t_j)+\bar V^j.
\end{equation}
The algorithm uses $2n$ noisy evaluations of $f$. Moreover, its combinatorial cost is $O(n)$ arithmetic operations.

In the case of exact information (i.e. $\delta=0$) we write $V^j$, $V^j_{\tau}$, $l$, $l_j$ instead of $\bar V^j$, $\bar V^j_{\tau}$, $\bar l$, $\bar l_j$, respectively. Of course $\bar l=\bar l_n$ but we will omit the subscript $n$ in order to simplify the notation.

Note that the algorithm \eqref{def_rand_RK_1} can be written as
\begin{equation}
\label{def_rand_RK_2}
		\left\{ \begin{array}{ll}
			\bar V^0 := \tilde \eta,\\
			\bar V_\tau^j := \bar V^{j-1} + h \tau_j \left[f (t_{j-1},  \bar V^{j-1}) +\delta_\tau^j \right], \\
			\bar V^j := \bar V^{j-1} + h\left[f(\theta_j, \bar V_\tau^j) +\delta^j \right],
		\end{array}\right.
\end{equation}
for $j\in\left\{1,\ldots,n\right\}$, where for the noise $\delta^j_{\tau}$, $\delta^j$ we have
\begin{equation}
\label{disturb_d1}
    \max_{1\leq j \leq n}\max\Bigl\{ \|\delta_\tau^j\|, \|\delta^j\|\Bigr\} \leq \delta
\end{equation}
almost surely. We stress that the only source of randomness in the noise are $\tau_j$'s, since $\delta^j$ is $\sigma(\tau_1,\ldots,\tau_j)$-measurable for $j\in\left\{1,\ldots,n\right\}$, while $\delta_{\tau}^j$ is $\sigma(\tau_1,\ldots,\tau_{j-1})$-measurable for $j\in\left\{2,\ldots,n\right\}$ ($\delta_{\tau}^1$ is deterministic). 

Lemma below is a crucial result that allow us to estimate error of the randomized Runge-Kutta algorithm in the case when $f$ is only locally Lipschitz. Similar localization technique was used in \cite{KaPr19}  for right-hand side functions that are globally Lipschitz but only locally differentiable.
\begin{lemma}
\label{Radii_R}
Let $\bar R=\bar R(a,b,K)$ be as in \eqref{def_bar_R}. Then for any  $n\in\mathbb{N}$, $\delta\in [0,1]$, $(\eta,f)\in F^{\varrho}$, and $(\tilde \eta,\tilde f)\in V_{(\eta,f)}(\delta)$  the following holds:
\begin{itemize}
    \item [(i)] There exists a unique solution $z=z(\eta,f)$, $z\in C^1([a,b];\mathbb{R}^d)$, of \eqref{PROBLEM}. Moreover, it holds that
    \begin{equation}
        z(t)\in B(\eta,\bar R), \quad t\in [a,b],
    \end{equation}
    there exists $\bar C=\bar C(a,b,K,L)\in (0,\infty)$ such that for all $s,t\in [a,b]$ 
    \begin{equation}
    \label{diff_z_2}
        \|z'(t)-z'(s)\|\leq \bar C |t-s|^{\varrho},
    \end{equation}
    and for all $j\in\left\{1,\ldots,n\right\}$ 
    \begin{equation}
        z(t_{j-1})+h\tau_jf(t_{j-1},z(t_{j-1}))\in B(\eta,\bar R)
    \end{equation}
    almost surely.
    \item[(ii)] For all $j\in\left\{0,1,\ldots,n\right\}$
    \begin{equation}
        V^j, \bar V^j\in B(\eta,\bar R)
    \end{equation}
    almost surely.
    \item[(iii)]For all $j\in\left\{1,\ldots,n\right\}$
    \begin{equation}
        V^j_{\tau}, \bar V^j_{\tau}\in B(\eta,\bar R)
    \end{equation}
    almost surely.
\end{itemize}
\end{lemma}
\noindent
{\bf Proof.} Let $(\eta,f)\in F^{\varrho}$. According to suitable version of Peano's theorem (see, for example, Theorem 70.4, page 292. in \cite{gorn_ing}) under conditions $(A1)$, $(A2)$ the equation \eqref{PROBLEM}  has at least one solution $z\in C^1([a,b];\mathbb{R}^d)$. Hence, firstly we show that there exists $R\in [0,\bar R]$  such that for any solution $z$ we have $ z(t) \in B(\eta,R)$ for all $t\in[a,b]$.
Note that
\begin{displaymath}
	\|z(t)\| \leq K(1 + b-a) + K \int\limits_a^t \|z(s)\|\rd s, \ t\in [a,b],
\end{displaymath}
and, by the Gronwall's lemma,
\begin{equation}
\label{est_sup_z_1}
    \sup\limits_{a\leq t\leq b}\|z(t) \| \leq C_1,
\end{equation}
where
$$ C_1 = C_1(a,b,K) = e^{K(b-a)} K(1+b-a). $$
Moreover, 
\begin{displaymath}
 	\sup\limits_{a\leq t\leq b}\|z(t) - \eta \| \leq  C_2,
\end{displaymath}
with
$$ C_2 = C_2(a,b,K) = C_1+K\leq \bar R,$$
where $\bar R$ is defined in \eqref{def_bar_R}. Hence, for all $z$ being the solution of the problem $\eqref{PROBLEM}$ and it holds
\begin{equation}
\label{z_incl_B}
    z(t) \in B(\eta,C_2)\subset B(\eta,\bar R), \ t \in [a,b].
\end{equation}
Now, let $z$ and $\tilde z$ are two solutions of \eqref{PROBLEM}. Due to \eqref{z_incl_B} we have that $z(t),\tilde z(t)\in  B(\eta,\bar R)$ for all $t\in [a,b]$. Therefore, by (A4) we have for all $t\in[a,b]$
$$ \|z(t)-\tilde z(t)\| \leq \int\limits_a^t \|f(s,z(s)) - f(s,\tilde z(s))\|\rd s \leq L \int\limits_a^t \|z(s)-\tilde z(s)\| \rd s.$$
This implies that for all $t\in[a,b]$ we have $z(t) = \tilde z(t)$, and the uniqueness in (i) follows. For the unique solution $z$ of \eqref{PROBLEM},  by $(A2)$ and \eqref{est_sup_z_1}, we have for all $t,s\in [a,b]$
\begin{displaymath}
    \|z(t)-z(s)\|\leq\int\limits_{\min\{t,s\}}^{\max\{t,s\}}\|f(u,z(u))\|du
    \leq K(1+C_1)|t-s|.
\end{displaymath}
Hence, by $(A3)$ and $(A4)$ we have for all $t,s\in [a,b]$
\begin{equation}
    \| z'(t) - z'(s)\| \leq L |t-s|^\rho + L \|z(t)-z(s)\|
    \leq \bar C|t-s|^{\varrho}
\end{equation}
 with $\displaystyle{\bar C= L\Bigl(1+K(1+C_1)(1+b-a)\Bigr)}$. This ends the proof of \eqref{diff_z_2}. Moreover, for all $n\in\mathbb{N}$, $j\in\left\{1,\ldots,n\right\}$, and almost surely
\begin{align*}
	\|z(t_{j-1}) + h\tau_j f(t_{j-1}, z(t_{j-1}))-\eta\| & \leq \|z(t_{j-1})-\eta\| + hK(1+\|z(t_{j-1})\|)\leq C_3,
\end{align*}
where $$C_3=C_3(a,b,K)=C_2 + (b-a) K(1+C_1).$$

For any $n\in\mathbb{N}$, $\delta\in [0,1]$, $(\tilde \eta,\tilde f)\in V_{(\eta,f)}(\delta)$, and $j\in\left\{0,1,\ldots,n\right\}$ we get
$$ \| \bar V^0\| \leq  K+1,
$$
\begin{equation}
\label{V_tau_est_1}
    \| \bar V_\tau^j \|  \leq \| \bar V^{j-1}\| + h \left( \|f(t_{j-1},\bar V^{j-1})\| + \|\delta_\tau^j\|\right)  \leq (1+hK)\| \bar V^{j-1}\| + h (K + 1),
\end{equation}
and therefore
\begin{displaymath}
    \|\bar V^j \| \leq \|\bar V^{j-1}\| + h\left( \|f(\theta_j, \bar V_\tau^j)\| +1 \right) \\
     \leq (1+K_1h) \|\bar V^{j-1}\| + K_2 h,
\end{displaymath}
where
\begin{displaymath}
 K_1 = K_1(a,b,K) := K(1+(b-a)K), \ K_2 = K_2(a,b,K) := (K+1)(1+(b-a)K).
\end{displaymath}
Hence, we get for all $0\leq j \leq n$ 
\begin{equation} 
\| \bar V^j \|  \leq (1+K_1h)^j (K+1) + \frac{(1+K_1h)^{j}-1}{K_1 h} K_2 h \leq C_4,
\end{equation}
where
\begin{displaymath}
     C_4=C_4(a,b,K):=\frac{K+1}{K}\Bigl((K+1)e^{K_1(b-a)}-1\Bigr).
\end{displaymath}
By \eqref{V_tau_est_1} we obtain for all $1\leq j \leq n$ 
\begin{equation} 
\| \bar V_\tau ^j \|  \leq C_5,
\end{equation}
with
\begin{displaymath}
     C_5=C_5(a,b,K):=C_4(1+(b-a)K) + (b-a)(K+1).
\end{displaymath}
Therefore
\begin{equation}
    \max\limits_{0\leq j\leq n}\|\bar V^j-\eta\|\leq \max\limits_{0\leq j\leq n}\|\bar V^j\|+\|\eta\|\leq C_6,
\end{equation}
\begin{equation}
    \max\limits_{1\leq j\leq n}\|\bar V_{\tau}^j-\eta\|\leq C_7,
\end{equation}
where
\begin{equation}
C_6=C_6(a,b,K)=C_4+K, \ C_7=C_7(a,b,K)=C_5+K.
\end{equation}
Note that $\bar R =\max\{C_2,C_3,C_6, C_7\}= \max \{C_3, C_7\}$ and the inclusions in (i)-(iii) follow. \ \ \ $\blacksquare$ \\ \\
We are ready to prove the following theorem that states upper error bounds for randomized Runge-Kutta algorithm under noisy information about the right-hand side function.
\begin{theorem} 
\label{error_RRK} Let $p\in [2,\infty)$. There exists $C$, depending only on the parameters of the class $F^{\varrho}$ and $p$, such that for all $n\geq \lfloor b-a\rfloor+1$, $\delta\in [0,1]$, $(\eta,f)\in F^{\varrho}$, $(\tilde \eta,\tilde f)\in V_{(\eta,f)}(\delta)$ we have
\begin{equation}
    \Bigl\|\sup\limits_{a\leq t\leq b}\|z(\eta,f)(t)-\bar l(\tilde \eta,\tilde f,\delta)(t)\|\Bigl\|_p\leq C\Bigl(h^{\varrho+1/2}+\delta\Bigr).
\end{equation}
\end{theorem}
\noindent
{\bf Proof.} We define
\begin{equation}
\label{def_bzj}
    \bar z_j(t)=\frac{z(t_{j+1})-z(t_j)}{h}(t-t_j)+z(t_j),
\end{equation}
for $t\in [t_j,t_{j+1}]$, $j\in\left\{0,\ldots,n-1\right\}$. Then
\begin{eqnarray}
\label{err_est_1}
    &&\Bigl\|\sup\limits_{a\leq t\leq b}\|z(t)-\bar l(t)\|\Bigl\|_p\leq \max\limits_{0\leq j\leq n-1}\sup\limits_{t_j\leq t\leq t_{j+1}}\|z(t)-\bar z_j(t)\|\notag\\
    &&\quad\quad +\Bigl\|\max\limits_{0\leq j\leq n-1}\sup\limits_{t_j\leq t\leq t_{j+1}}\|\bar z_j(t)-\bar l_j(t)\|\Bigl\|_p,
\end{eqnarray}
where, from \eqref{def_blj}, \eqref{def_bzj}, it holds
\begin{eqnarray}
&&\Bigl\|\max\limits_{0\leq j\leq n-1}\sup\limits_{t_j\leq t\leq t_{j+1}}\| \bar z_j(t) - \bar l_j(t) \|\Bigl\|_p \leq 3\Bigl\|\max\limits_{0\leq j\leq n} \|z(t_j) - V^j\|\Bigl\|_p\notag\\
&&\quad\quad +3\Bigl\|\max\limits_{0\leq j\leq n} \|V^j - \bar V^j\|\Bigl\|_p.
\end{eqnarray}
In the case of global Lipschitz condition the first term  $\Bigl\|\max\limits_{0\leq j\leq n} \|z(t_j) - V^j\|\Bigl\|_p$ was estimated in Theorem 5.2 in \cite{KruseWu_1}. Under local Lipschitz assumptions $(A3),(A4)$ and with the help of Lemma \ref{Radii_R} we can estimate it essentially in the same fashion as in \cite{KruseWu_1}, and we get the same upper bound
\begin{equation}
\label{upp_ex_inf}
    \Bigl\|\max\limits_{0\leq j\leq n} \|z(t_j) - V^j\|\Bigl\|_p\leq Ch^{\varrho+1/2}.
\end{equation}
 However, for the convenience of the reader we present complete justification of \eqref{upp_ex_inf}, where we explicitly point out the use of Lemma \ref{Radii_R}.
 
For $k\in\left\{1,\ldots,n\right\}$ we have
\begin{equation}
\label{cr_1}
    z(t_k)-V^k=S_1^k+S_2^k+S_3^k,
\end{equation}
where
\begin{equation}
    S_1^k=\sum\limits_{j=1}^k\Bigl(\int\limits_{t_{j-1}}^{t_j} z'(s) ds-h z'(\theta_j)\Bigr),
\end{equation}
\begin{equation}
    S_2^k=h\sum\limits_{j=1}^k\Bigl(f(\theta_j,z(\theta_j))-f(\theta_j,z(t_{j-1})+h\tau_jf(t_{j-1},z(t_{j-1})))\Bigr),
\end{equation}
\begin{equation}
    S_3^k=h\sum\limits_{j=1}^k\Bigl(f(\theta_j,z(t_{j-1})+h\tau_j f(t_{j-1},z(t_{j-1})))-f(\theta_j,V^j_{\tau})\Bigr).
\end{equation}
By \eqref{diff_z_2} in Lemma \ref{Radii_R} and Theorem 3.1 in \cite{KruseWu_1} we have that there exists $C_1>0$ such that for all $n\geq \lfloor b-a\rfloor+1$
\begin{equation}
\label{cr_2}
    \Bigl\|\max\limits_{1\leq k\leq n}\|S_1^k\| \Bigl\|_p=\Biggl\|\max\limits_{1\leq k\leq n}\Bigl\|\int\limits_{a}^{t_k}z'(s)ds-h\sum\limits_{j=1}^k z'(\theta_j)\Bigl\| \Biggl\|_p    \leq C_1h^{\varrho+1/2}.
\end{equation}
Furthermore, by Lemma \ref{Radii_R} and $(A4)$ we get for $k\in\{1,\ldots,n\}$
\begin{eqnarray}
\label{cr_3}
    &&\|S_2^k\|\leq hL\sum\limits_{j=1}^k\|z(\theta_j)-z(t_{j-1})-h\tau_jf(t_{j-1},z(t_{j-1}))\|\notag\\
    &&\leq hL\sum\limits_{j=1}^k\int\limits_{t_{j-1}}^{\theta_j}\|z'(s)-z'(t_{j-1})\|ds\leq C_2h^{\varrho+1}.
\end{eqnarray}
Moreover, by \eqref{def_rand_RK_2} (with $\delta=0$), Lemma \ref{Radii_R} and $(A4)$ we have for $k\in\{1,\ldots,n\}$
\begin{eqnarray}
\label{cr_4}
    &&\|S_3^k\|\leq hL\sum\limits_{j=1}^k\|z(t_{j-1})+h\tau_j f(t_{j-1},z(t_{j-1}))-V^j_{\tau}\|\notag\\
    &&\leq hL(1+L(b-a))\sum\limits_{j=1}^k\|z(t_{j-1})-V^{j-1}\|\leq hC_3\sum\limits_{j=1}^k\max\limits_{0\leq i\leq j-1}\|z(t_{i})-V^{i}\|.
\end{eqnarray}
From \eqref{cr_1}, \eqref{cr_2}, \eqref{cr_3}, and \eqref{cr_4} we have for $k\in\{1,\ldots,n\}$ that
\begin{eqnarray}
    &&\Bigl\|\max\limits_{1\leq i\leq k}\|z(t_i)-V^i\|\Bigl\|_p\leq \Bigl\|\max\limits_{1\leq i\leq n}\|S_1^i\|\Bigl\|_p+\Bigl\|\max\limits_{1\leq i\leq n}\|S_2^i\|\Bigl\|_p+hC_3\sum\limits_{j=1}^{k-1}\Bigl\|\max\limits_{0\leq i\leq j}\|z(t_i)-V^i\|\Bigl\|_p\notag\\
    &&\leq C_4h^{\varrho+1/2}+hC_3\sum\limits_{j=1}^{k-1}\Bigl\|\max\limits_{0\leq i\leq j}\|z(t_i)-V^i\|\Bigl\|_p.
\end{eqnarray}
Using the weighted version of discrete Gronwall's lemma (see, for example, Lemma 2.1. in \cite{KruseWu_1}) we obtain \eqref{upp_ex_inf}.

We now establish upper bound on $\Bigl\|\max\limits_{0\leq j\leq n} \|V^j - \bar V^j\|\Bigl\|_p$. We have 
$$ \| V^0 - \bar V^0\|  \leq \delta, $$
and by Lemma \ref{Radii_R} (ii), (iii) it holds for $j\in\left\{1,\ldots,n\right\}$
\begin{equation}
    \|V^j - \bar V^j\|\leq \|V^{j-1} - \bar V^{j-1}\| + h L \|V_\tau^j - \bar V_\tau^j\| + h\delta,
\end{equation}
where
\begin{equation}
    \| V_\tau^j - \bar V_\tau^j \| \leq 
     \|V^{j-1} - \bar V^{j-1}\| + h\delta + h L \|V^{j-1} - \bar V^{j-1}\|.
\end{equation}
Hence
\begin{equation}
    \| V^j - \bar V^j \| \leq (1+C_5h)\cdot \|V^{j-1} - \bar V^{j-1}\| + h\delta C_6,
\end{equation}
and for all $j\in\left\{0,1,\ldots,n\right\}$
\begin{equation}
    \|V^j - \bar V^j\|  \leq (1+C_5h)^j \|V^0 - \bar V^0\| + \frac{(1+C_5h)^j-1}{C_5 h}\cdot h\delta C_6\leq C_7\delta.
\end{equation}
Therefore,
\begin{equation}
\label{zero_stab}
    \max_{0\leq j \leq n} \|V^j-\bar V^j\| \leq C_7 \delta,
\end{equation}
almost surely and by \eqref{upp_ex_inf} we get
\begin{equation} 
\label{est_zjblj}
   \Bigl\|\max\limits_{0\leq j\leq n-1}\sup\limits_{t_j\leq t\leq t_{j+1}}\| \bar z_j(t) - \bar l_j(t) \|\Bigl\|_p  \leq C_8(h^{\varrho+1/2} + \delta).
\end{equation}
We now show the upper bound on $\max\limits_{0\leq j\leq n-1}\sup\limits_{t_j\leq t\leq t_{j+1}}\|z(t)-\bar z_j(t)\|$, which is the deterministic term in the error estimate \eqref{err_est_1}.

For every $t\in[t_j, t_{j+1}]$, $j\in\left\{0,\ldots,n-1\right\}$ we obtain,  by applying the mean value-theorem  component-wise,
\begin{displaymath}
z(t) = \sum_{k=1}^d z_k(t)e_k = \sum_{k=1}^d\Bigl (z_k(t_j)
+z_k'(\alpha_{k,j}^t)(t-t_j)\Bigr)e_k,
\end{displaymath}
\begin{equation}
    \bar z_j(t)=\sum_{k=1}^d \bar z_{j,k}(t)e_k = \sum_{k=1}^d \Bigl(z_k'(\beta_{k,j})(t-t_j) + z_k(t_j)\Bigr)e_k
\end{equation}
for some $\alpha^t_{k,j}\in[t_j,t]\subset[t_j, t_{j+1}]$ and $ \beta_{k,j} \in [t_j,t_{j+1}].$ Thereby, for $t\in [t_j,t_{j+1}]$
\begin{equation} 
z(t) - \bar z_j(t)  = \sum_{k=1}^d \Bigl( z_k'(\alpha^t_{k,j}) - z_k'(\beta_{k,j}) \Bigr)(t-t_j)e_k.
\end{equation}
Since $\|e_k\|=1$ for $k\in\left\{1,\ldots,d\right\}$, by \eqref{diff_z_2} we get for $j\in\left\{0,\ldots,n-1\right\}$,  $t\in [t_j,t_{j+1}]$
\begin{equation}
\label{est_zbzj}
 \|z(t)-\bar z_j(t) \| \leq  h\sum_{k=1}^d |z_k'(\alpha_{k,j}^t)-z_k'(\beta_{k,j})|\leq d\bar C h^{\varrho+1}.
\end{equation}
Combining \eqref{err_est_1}, \eqref{est_zjblj}, and \eqref{est_zbzj} we get the thesis. \ \ \ $\blacksquare$
\section{Lower bounds and optimality of the randomized Runge-Kutta algorithm}
This section is devoted to  lower bounds on the worst-case error of any algorithm from the class $\Phi_n$. They will allow us to conclude that the randomized Runge-Kutta $\bar l$ is asymptotically optimal within this setting. 

\begin{lemma}\label{lower-bound-lemma}
	Let $p\in [2,\infty)$ and $\varrho\in (0,1]$, then
	$$    e^{(p)}_n({F^{\varrho}},\delta)  = \Omega(\max\{n^{-(\varrho+1/2)},\delta\}) $$
as $n \to \infty$ and $\delta \to 0+$.
\end{lemma}
\noindent
{\bf Proof.}
    Firstly, for the exact randomized information the following lower bound holds 
    \begin{equation}
    \label{lb_1}
     e^{(p)}_n({F^{\varrho}}, \delta)\geq e^{(p)}_n({F^{\varrho}}, 0)  = \Omega(n^{-(\varrho+1/2)}), \ n \to \infty.
     \end{equation}
    This follows from reducing an integration problem of H\"older continuous functions to the solution of initial value problem, see \cite{HeinMilla} and \cite{NOV} for the details.
  
    Note that for any algorithm $\mathcal{A}\in\Phi_n$ and any $(\eta_1,f_1),(\eta_2,f_2)\in F^{\varrho}$, such that $V_{(\eta_1,f_1)}(\delta)\cap V_{(\eta_2,f_2)}(\delta)\neq\emptyset$, we have
    \begin{equation}
        e^{(p)}(\mathcal{A},F^{\varrho},\delta)\geq \frac{1}{2}\sup\limits_{a\leq t\leq b}\|z(\eta_1,f_1)(t)-z(\eta_2,f_2)(t)\|.
    \end{equation}
Hence, let us take  $(\eta_1,f_1)=(0e_1,+\delta e_1)$, $(\eta_2,f_2)=(0e_1,-\delta e_1)$ that belong to $F^{\varrho}$ if $\delta\in [0,\min\{K,1\}]$.  Then  $(0e_1,0e_1)\in V_{(\eta_1,f_1)}(\delta)\cap V_{(\eta_2,f_2)}(\delta)$  and
    $$
        e^{(p)}(\mathcal{A},F^{\varrho},\delta)\geq 
          \delta\sup\limits_{a\leq t\leq b}\|  (t-a)e_1 \| = (b-a)\delta,
    $$
    which implies the following
    \begin{equation}
    \label{lb_2}
        e^{(p)}_n({F^{\varrho}}, \delta)\geq (b-a)\delta.
    \end{equation}
    By \eqref{lb_1} and \eqref{lb_2} we get the thesis. \ \ \ $\blacksquare$ \\ \\
Lemma \ref{lower-bound-lemma}, together with Theorem \ref{error_RRK} immediately imply the following theorem of optimality of randomized Runge-Kutta algorithm.
\begin{theorem}
\label{thm_opt_RK}
    Let $p\in [2,\infty)$ and $\varrho\in (0,1]$, then
	$$    e^{(p)}_n({F^{\varrho}},\delta)  = \Theta(\max\{n^{-(\varrho+1/2)},\delta\}) $$
    as $n \to \infty$ and $\delta\to 0+$. The optimal algorithm is the randomized Runge-Kutta algorithm $\bar l$.
\end{theorem}
\begin{remark}
    If we restrict considerations to {\it deterministic algorithms}, then  the following sharp bounds on the $n$th minimal error hold
    \begin{equation}
         e^{(p)}_n({F^{\varrho}},\delta)  = \Theta(\max\{n^{-\varrho},\delta\}),
    \end{equation}
    as $n \to \infty$ and $\delta \to 0+$. The classical Euler scheme, based on the equidistant mesh,  is the optimal one.
\end{remark}
\section{Numerical experiments}
\begin{figure}[h]
\centering
\includegraphics[width=0.4\linewidth]{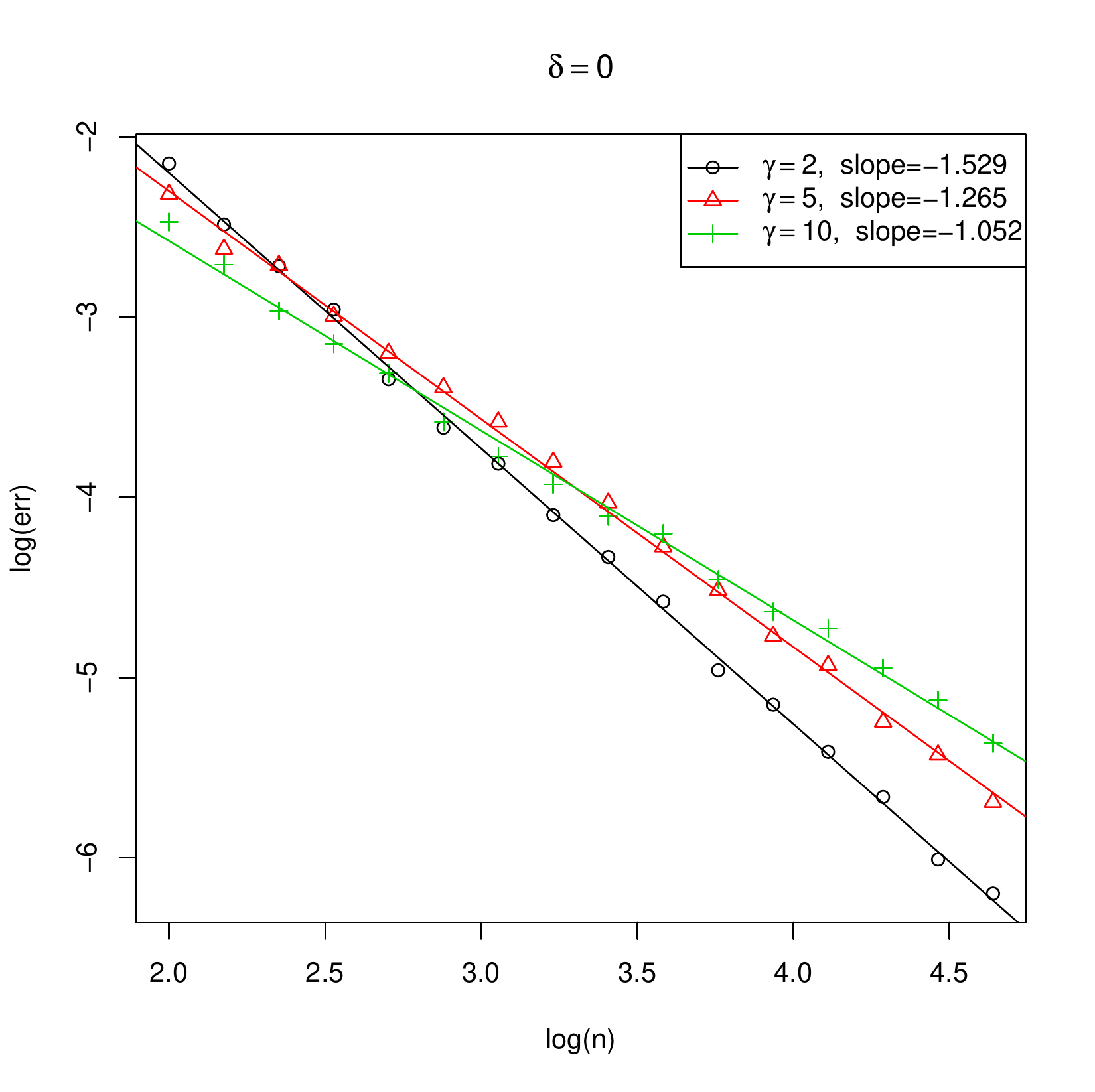}
\includegraphics[width=0.4\linewidth]{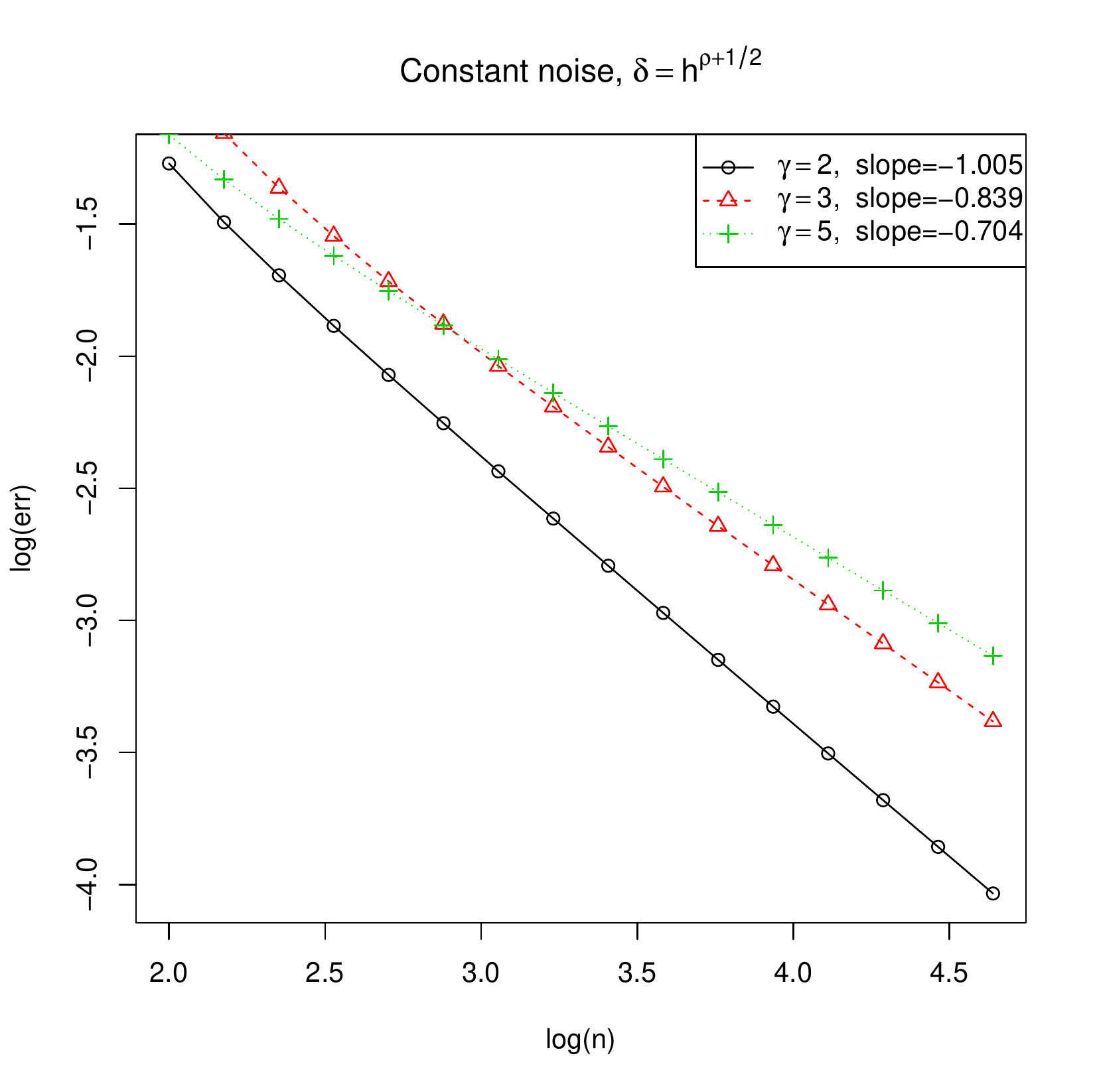}
\caption{$\lg(err)$ vs $\lg(n)$ for Example 1 with $\delta=0$ (left graph) and $\delta=h^{\rho+1/2}$ (right graph)}
\label{fig:problem1_delta0}\end{figure}

In order to support the obtained theoretical results we conducted several numerical experiments. The worst case noise was simulated in two ways: by using two constant noises equal to $\delta$ and $-\delta$ (as in the proof of the lower bounds), and then taking the worst of them, and by generating $100$ of repetitions of random noise from the uniform distribution on $[-\delta,\delta]$ for each step of the algorithm, and then by taking the worst of them. The error was approximated in $L^2$ norm at the terminal point by $M$ repetitions of the randomized Runge-Kutta scheme (with $M$ equal to $1000$ for the constant noise and $100$ for the random noise). 
\newline\newline\bf Example 1. \rm 
As the first problem we consider  the following scalar ODE
\begin{equation}
	\label{problem 1}
		\left\{ \begin{array}{ll}
			z'(t)= 1+z(t) \cos\left(10 (2-t) ^{1/\gamma}|z(t)|^{3/2}\right),\quad t\in[0,2]\\
			z(0) = -1 
		\end{array}\right.
\end{equation}
for different values of $\gamma\in (0,+\infty)$. Note that the right-hand side function $f$ in \eqref{problem 1} satisfies the assumptions (A1)--(A4) with $\rho=1/\gamma$.  
The results for exact information ($\delta=0$), $\gamma\in\left\{2,5,10\right\}$ and $n$ varying for 100 to 50000 are presented in Figure~\ref{fig:problem1_delta0} (left graph). The results are printed in the logarithmic scale as the relation $\lg(err)$ versus $\lg(n)$. We have added also the slope of the relations. Note that we get a little better behavior that it follows from the theoretical results, which may be due to the fact that the the right-hand side function $f=f(t,y)$ is Lipschitz continuous with respect to time variable $t$ on every interval $[0,\beta]$ with $\beta<2$. 

\begin{figure}[h]
\centering
\includegraphics[width=0.4\linewidth]{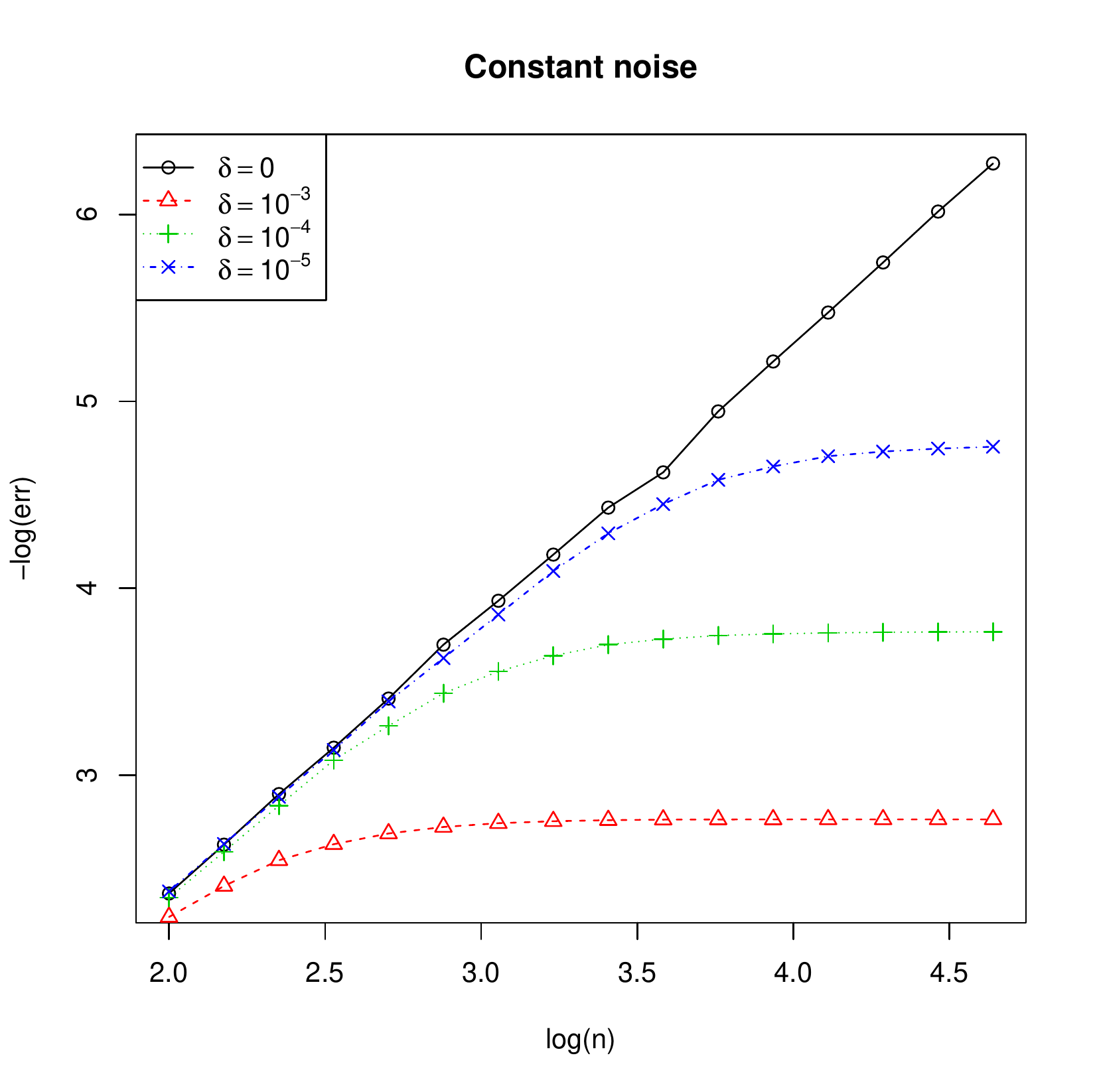}
\includegraphics[width=0.4\linewidth]{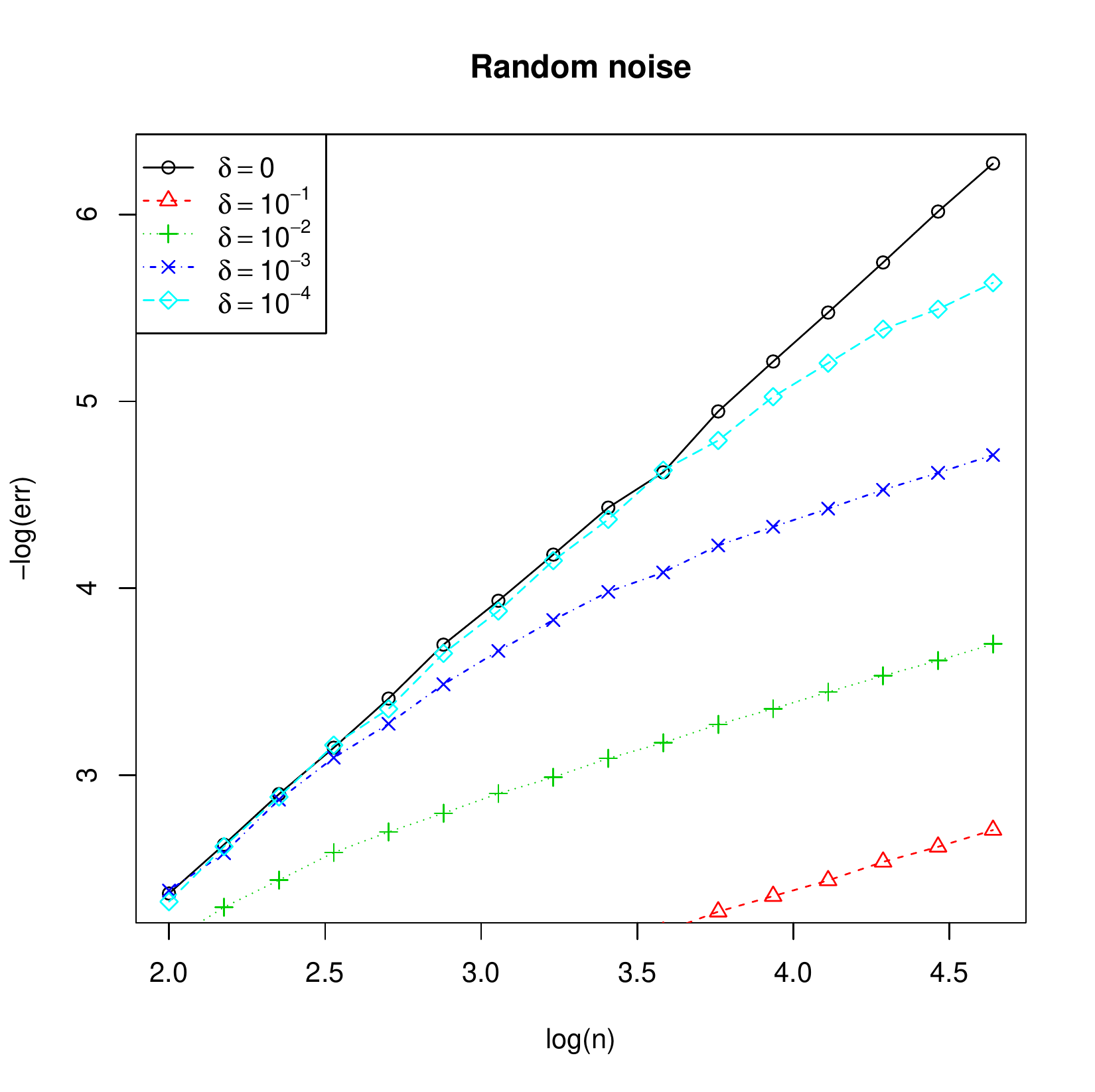}
\caption{$-\lg(err)$ vs $\lg(n)$ for Example 1, $\gamma=3$ (left graph -- constant noise, right graph -- random noise)}
\label{fig:problem1_gamma3}\end{figure}

In Figure~\ref{fig:problem1_gamma3} the relation $-\lg(err)$ versus $\lg(n)$ for $\gamma=3$ and different values of $\delta$ is presented.
We have also run the test with varying $ \delta \asymp h^{\rho +1/2}$. The result is presented in Figure~\ref{fig:problem1_delta0} (right graph). As we can see, the error decreases proportionally to $\delta$, which confirms the theoretical results.
\newline\newline\bf Example 2. \rm 
Below we recall the well-known  SIR model that models the spread of disease
\begin{equation}
	\label{problem 2}
			[S'(t), I'(t), R'(t)] = [- \beta S(t) I(t), \beta S(t) I(t) - \gamma I(t),\gamma I(t)], \ t\in [0,30].
\end{equation}
In numerical experiments we set $[S(0), I(0), R(0)]=[50,1,0]$, $\beta = \frac{1}{768}$, $\gamma = \frac{1}{120}$. The right-hand side function in \eqref{problem 2}  does not belong to the class $F^{\rho}$, since it is not globally of at most linear growth. However, we still achieve the desired empirical convergence rate $O(h^{3/2})$. This suggests possibility of  weakening the assumption (A2) in the future investigations.

\begin{figure}[h]
\centering
\includegraphics[width=0.4\linewidth]{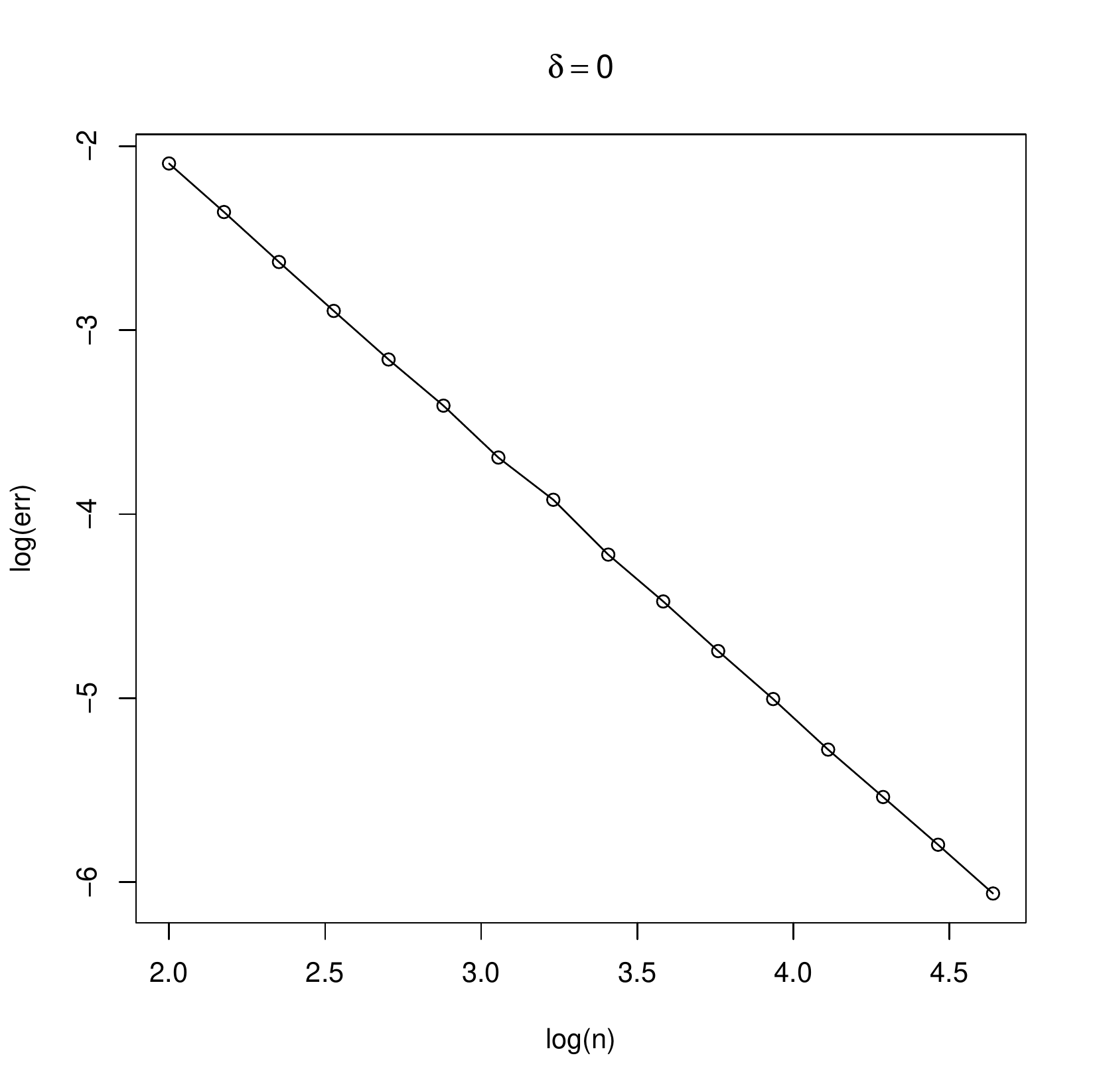}
\includegraphics[width=0.4\linewidth]{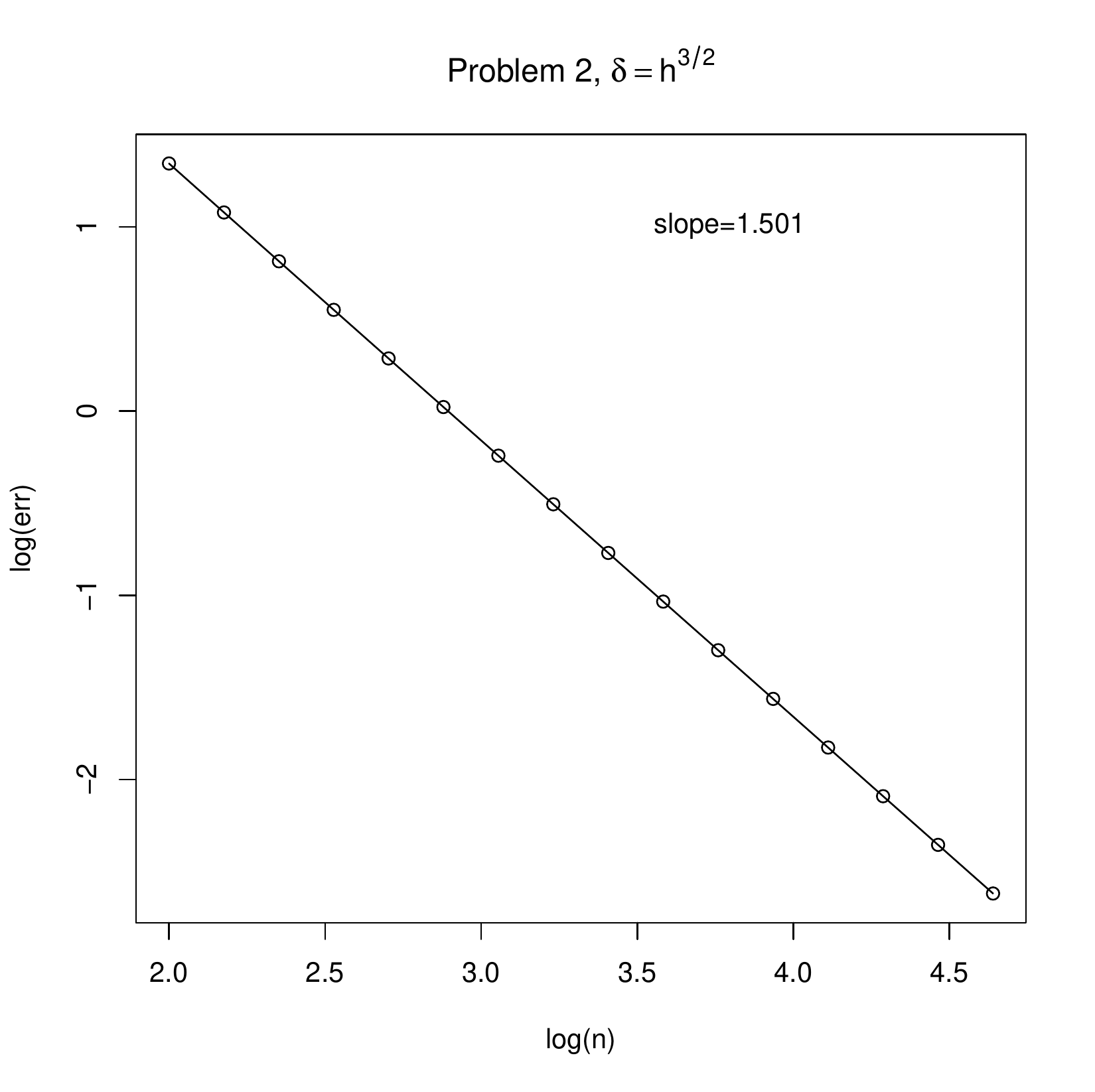}
\caption{$\lg(err)$ vs $\lg(n)$ for Example 2 with $\delta=0$ (left graph) and $\delta=h^{3/2}$ (right graph)}
\label{fig:problem2_delta0}\end{figure}

\begin{figure}[h]
\centering
\includegraphics[width=0.4\linewidth]{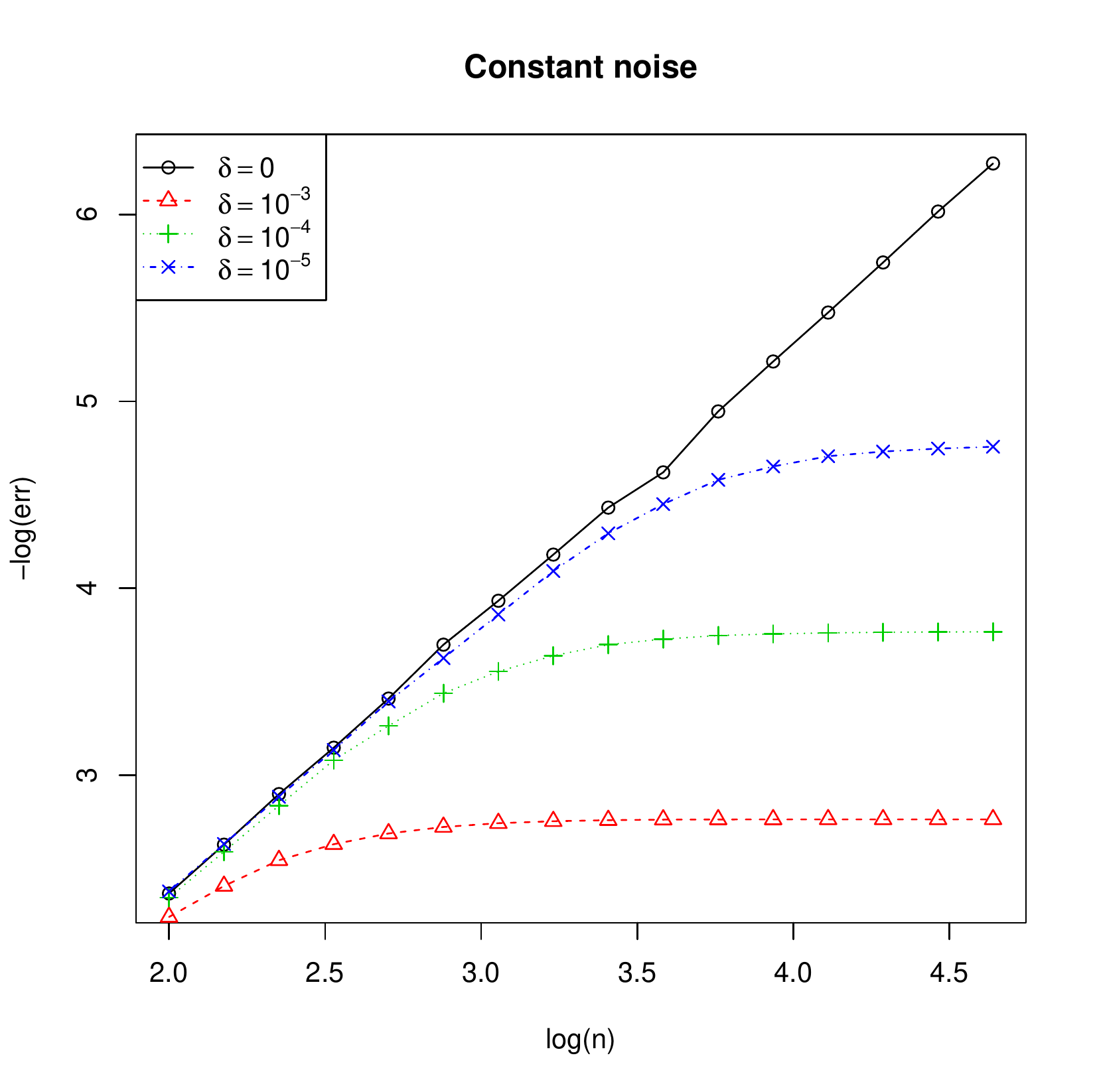}
\includegraphics[width=0.4\linewidth]{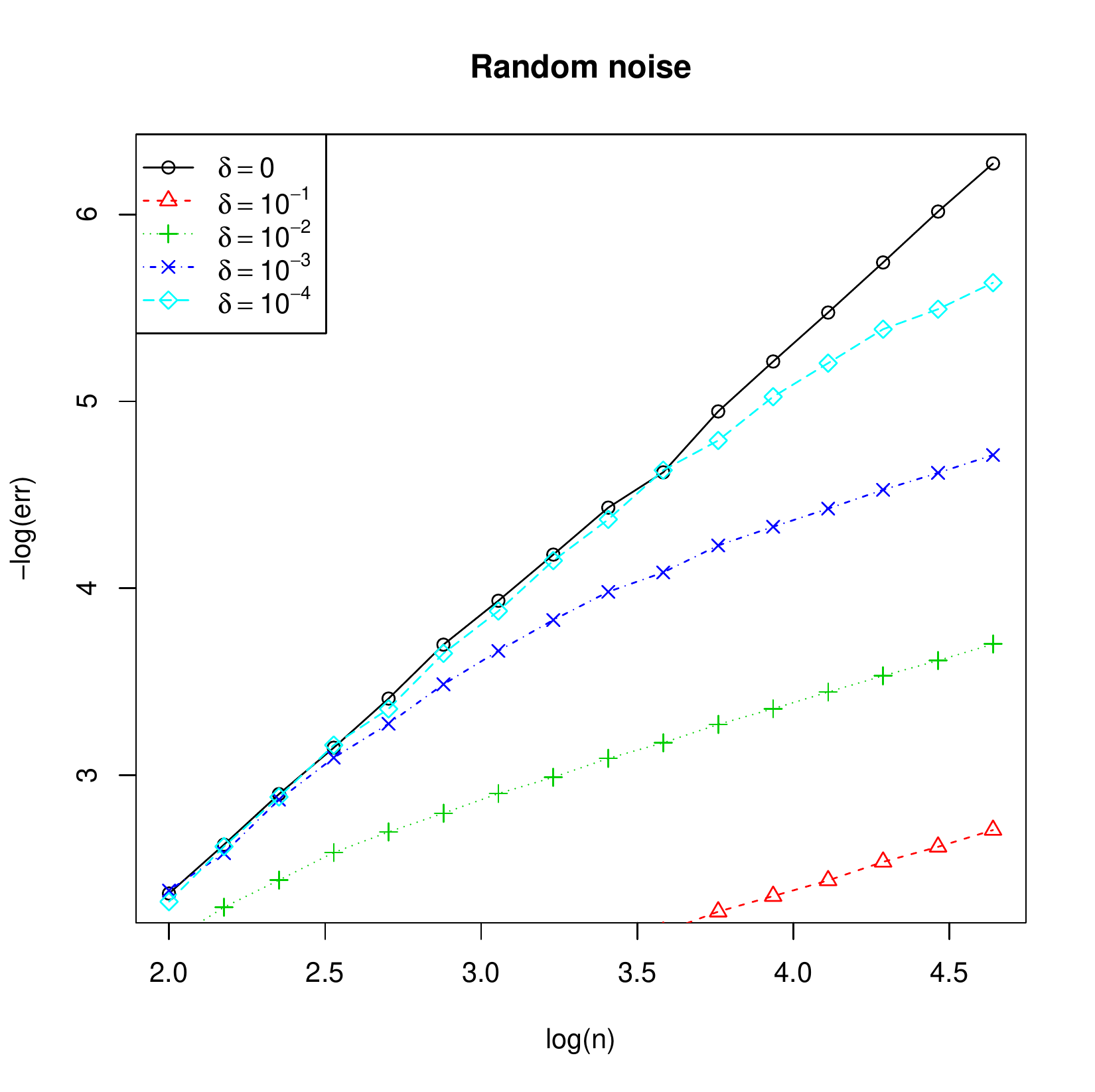}
\caption{$-\lg(err)$ vs $\lg(n)$ for Example 2 (left graph -- constant noise, right graph -- random noise)}
\label{fig:problem2}\end{figure}

We have made similar simulation as for Example 1. In Figure~\ref{fig:problem2_delta0} we present the relation $\lg(err)$ versus $\lg(n)$ for exact information (left graph) and for inexact information with the precision parameter $\delta\asymp h^{3/2}$ (right graph). We can see on both graphs that the error is proportional to $n^{-3/2}$. In Figure~\ref{fig:problem2} we present also the relation $-\lg(err)$ versus $\lg(n)$ for different values of $\delta$.

\section{Three types of regions of stability for the randomized Runge-Kutta method}\label{sec:stab}
In the case when the information is exact we investigate absolute stability of the randomized Runge-Kutta method $(V^j)_{j\in\left\{0,1,\ldots,n\right\}}$. Since now the algorithm is random we have to generalize definitions concerning the absolute stability known for deterministic methods, see, for example, \cite{hawan}.

Let us consider the well-known  test problem
\begin{equation}
	\label{TEST_PROBLEM}
		\left\{ \begin{array}{ll}
			z'(t)= \lambda z(t), \ t\geq 0, \\
			z(0) = \eta
		\end{array}\right.
\end{equation}
with $\lambda\in\mathbb{C}$, $\eta\neq 0$. The exact solution of \eqref{TEST_PROBLEM} is $z(t)=\eta\exp(\lambda t)$ and 
\begin{equation}
    \lim\limits_{t\to\infty} z(t)=0 \ \hbox{iff} \ \Re(\lambda)<0.
\end{equation}
For a fixed step-size $h>0$ we apply the algorithm $(V^j)_{\in\left\{0,1,\ldots,n\right\}}$ based on the mesh $t_j=jh$, $j\in\mathbb{N}_0$, to the test problem \eqref{TEST_PROBLEM}. As a result we obtain the following recurrence
\begin{equation}
\label{REC_1}
    V^0 = \eta \ \text{and} \ V^j=p_j(h\lambda)\cdot V^{j-1}  \ \text{for} \ j\in\left\{1,\ldots,n\right\},
\end{equation}
where 
\begin{equation}
\label{def_p_j}
    p_j(z) = \tau_j\cdot z^2+z+1, \ z\in\mathbb{C}
\end{equation}
    is a second-degree polynomial with random coefficient $\tau_j$. Let us substitute $z=h\lambda$. For any $z\in\mathbb{C}$, $(p_j(z))_{j\in\mathbb{N}}$ is a sequence of complex-valued, independent, and identically distributed random variables on $(\Omega, \Sigma, \mathbb{P})$. Solving \eqref{REC_1} we get that
\begin{equation}
\label{REC_2}
    V^k = \eta\cdot\prod\limits_{j=1}^k p_j(z).
\end{equation}
We consider three sets
\begin{eqnarray}
    &&\mathcal{R}^{MS} = \{z\in\mathbb{C} \colon V^k\to 0 \ \hbox{in $L^2(\Omega)$ as} \ k\to \infty \}, \notag\\
    &&\mathcal{R}^{AS} = \{z\in\mathbb{C} \colon V^k\to 0 \  \hbox{almost surely as} \ k\to \infty\},\notag\\
     &&\mathcal{R}^{SP} = \left\{ z\in \mathbb{C} \colon  V^k \to 0 \text{ in probability as} \ k\to \infty  \right\},
    \label{eq:regs}
\end{eqnarray}
where we call $\mathcal{R}^{MS}$ the region of mean-square stability,  $\mathcal{R}^{AS}$ -- the region of asymptotic stability, while $\mathcal{R}^{SP}$ -- the region of stability in probability. Of course we have that 
\begin{equation} 
\label{reg_incl}
    \mathcal{R}^{MS}\cup\mathcal{R}^{AS}\subset\mathcal{R}^{SP},
\end{equation} but we will show more accurate inclusions. If in $(V^j)_{j\in\left\{0,1,\ldots,n\right\}}$ we set  $\tau_j:=1/2$ for all $j\in\mathbb{N}$ then we arrive at the well-known (deterministic) midpoint scheme. The well-known region of absolute stability of this algorithm is 
\begin{equation}
    \mathcal{R}^{Mid}=\Bigl\{z\in\mathbb{C} \colon \Bigl|\frac{1}{2}z^2+z+1\Bigr|^2<1\Bigr\}.
\end{equation}
We use $\mathcal{R}^{Mid}$ as a reference set for $\mathcal{R}^{MS}$, $\mathcal{R}^{AS}$, and $\mathcal{R}^{SP}$. We also investigate the intervals of absolute stability
\begin{equation}
    \mathcal{I}^{\diamond}=\mathcal{R}^{\diamond}\cap\{z\in\mathbb{C} \colon  \Im(z)=0\}, \ \diamond\in\{MS,AS,SP\}.
\end{equation}
In \cite{higham2000}  regions $\mathcal{R}^{MS}$, $\mathcal{R}^{AS}$ were defined in order to investigate stability of numerical schemes for stochastic differential equations. Here we adopt this methodology for randomized Runge-Kutta algorithm in the context of deterministic ordinary differential equations. According to our best knowledge this is the first attempt in that direction. Moreover, we also investigate properties of the region $\mathcal{R}^{SP}$, which was not the case in \cite{higham2000}.
\subsection{Region $\mathcal{R}^{MS}$ of mean-square stability}
For all $z\in\mathbb{C}$ and $j\in\mathbb{N}$
\begin{equation}
\label{pj_sq}
    |p_j(z)|^2=1+2\Re(z)+|z|^2+2\tau_j \Bigl(2(\Re(z))^2-|z|^2+|z|^2 \Re(z)\Bigr)+\tau_j^2|z|^4,
\end{equation}    
and hence
\begin{equation}
\label{exp_pj}
    \mathbb{E}|p_j(z)|^2=1+2\Re(z)\Bigl(1+\frac{1}{2}|z|^2\Bigr)+2(\Re(z))^2+\frac{1}{3}|z|^4.
\end{equation}   
By \eqref{REC_2} we get that
\begin{equation}
    \|V^k\|_2=|\eta|\cdot (\mathbb{E}|p_1(z)|^2)^{k/2}.
\end{equation}
Hence, we can write that
\begin{equation}
\label{def_MS_reg}
    \mathcal{R}^{MS}=\{z\in\mathbb{C} \colon \mathbb{E}|p_1(z)|^2<1\}.
\end{equation}
\begin{theorem}\noindent
\label{prop_MS}
    \begin{itemize}
        \item [(i)]  The sets $\mathcal{R}^{MS}$, $\mathcal{R}^{Mid}$ are open and  symmetric with respect to the real axis.
        \item [(ii)] There exists $r_0\in (0,\infty)$ such that $\mathcal{R}^{MS}\subset\mathcal{R}^{Mid}\subset \mathbb{C}_{-}\cap \{z\in \mathbb{C} \colon |z|<r_0\}$.
        \item [(iii)] $\mathcal{I}^{MS}=(x_0,0)$ with $\displaystyle{x_0=-1-(\sqrt{2}-1)^{-1/3}+(\sqrt{2}-1)^{1/3}}$, while $\mathcal{I}^{Mid}:=\mathcal{R}^{Mid}\cap\{z\in\mathbb{C} \colon  \Im(z)=0\}=(-2,0)$, and $\displaystyle{\mathcal{I}^{MS}\subset \mathcal{I}^{Mid}}$.
    \end{itemize}
\end{theorem}
\noindent
{\bf Proof.} It holds 
\begin{eqnarray}
    &&\mathcal{R}^{MS}=(\phi^{MS})^{-1}((-\infty,0)),\\
    &&\mathcal{R}^{Mid}=(\phi^{Mid})^{-1}((-\infty,0)),
\end{eqnarray}
where the functions $\phi^{MS},\phi^{Mid}:\mathbb{C}\to\mathbb{R}$ are given as follows
\begin{eqnarray}
    &&\phi^{MS}(z)=2\Re(z)\Bigl(1+\frac{1}{2}|z|^2\Bigr)+2(\Re(z))^2+\frac{1}{3}|z|^4,\notag\\
    &&\phi^{Mid}(z)=2\Re(z)\Bigl(1+\frac{1}{2}|z|^2\Bigr)+2(\Re(z))^2+\frac{1}{4}|z|^4.
\end{eqnarray}
Note that $\phi^{MS},\phi^{Mid}$ are continuous, thereby  $\mathcal{R}^{MS}$, $\mathcal{R}^{Mid}$ are open.
Since for all $z\in\mathbb{C}$ we have that $\phi^{MS}(z)=\phi^{MS}(\bar z)$, $\phi^{Mid}(z)=\phi^{Mid}(\bar z)$, the conclusion in (i) follows. Moreover $\phi^{Mid}(z)\leq \phi^{MS}(z)$ for all $z\in\mathbb{C}$, and hence $\mathcal{R}^{MS}\subset\mathcal{R}^{Mid}$. Furthermore, for any $z\in\mathcal{R}^{Mid}$  we have
\begin{equation}
    \Re(z)<\frac{-2(\Re(z))^2-\frac{1}{4}|z|^4}{2\Bigl(1+\frac{1}{2}|z|^2\Bigr)}<0,
\end{equation}
which implies that
\begin{equation}
\label{R_MID_INCL_1}
    \mathcal{R}^{Mid}\subset\mathbb{C}_{-}.
\end{equation}
Let us consider any $z\in\mathcal{R}^{Mid}$. Then, by the definition of $\mathcal{R}^{Mid}$ we obtain what follows
\begin{equation}
    1 > \left| \frac12 z^2+z+1 \right| \geq \left| \frac12 z^2+z \right| - 1 = |z|\cdot \left| \frac12 z + 1 \right| - 1 \geq |z|\cdot \left( \frac12 \left|  z \right|- 1\right)  - 1,
\end{equation}
leading to 
\begin{equation}
    0 > \frac12 |z|^2 - |z| -2 = \frac12 \cdot \left( |z| - \left( 1+\sqrt{5} \right) \right)\cdot \left( |z| - \left( 1-\sqrt{5} \right) \right).
\end{equation}
Since $|z| - \left( 1-\sqrt{5} \right)>0$ for all $z\in\mathbb{C}$, we conclude that $|z|<1 + \sqrt{5}$. Thus, 
\begin{equation} \label{R_MID_INCL_2}
    \mathcal{R}^{Mid} \subset \left\{ z\in \mathbb{C} \colon |z|<r_0 \right\},
\end{equation}
where $r_0 = 1 + \sqrt{5}$. Inclusion \eqref{R_MID_INCL_1}, combined with \eqref{R_MID_INCL_2}, leads to (ii).
By \eqref{def_MS_reg} and Cardano's formula we get that
\begin{equation}
    \mathcal{I}^{MS}=\{x\in (-\infty,0) \colon x^3+3x^2+6x+6>0\}=(x_0,0),
\end{equation}
where
\begin{equation}
    x_0=-1-\frac{1}{(\sqrt{2}-1)^{1/3}}+(\sqrt{2}-1)^{1/3}.
\end{equation}
It is well-known that for the deterministic mid-point method $\mathcal{I}^{Mid}=(-2,0)$, however, for the convenience of the reader we provide a short justification. Namely, from \eqref{pj_sq} with $\tau_j=1/2$ we get that
\begin{equation}
    \mathcal{I}^{Mid}=\{x\in (-\infty,0) \colon x^3+4x^2+8x+8>0\}=(-2,0).
\end{equation}
Since $x_0>-2$ we get the inclusion $\mathcal{I}^{MS}\subset\mathcal{I}^{Mid}$. This proves (iii). \ \ \ $\blacksquare$ 
\subsection{Region $\mathcal{R}^{AS}$ of asymptotic stability}
The following result is a rearrangement of Lemma 5.1 in \cite{higham2000}.
\begin{lemma} 
\label{as_stab_lem_1}
    Given a sequence of real-valued, independent and identically distributed random variables $\{Z_n\}_{n\in\mathbb{N}_0}$ with $\displaystyle{\mathbb{P}(Z_1>0)=1}$, consider the sequence of random variables $\{Y_n\}_{n\in\mathbb{N}}$ defined by
    \begin{equation}
        Y_n=\Bigl(\prod\limits_{i=0}^{n-1}Z_i\Bigr)Y_0,
    \end{equation}
    where $Y_0$ is independent of $\{Z_n\}_{n\in\mathbb{N}_0}$ and $\displaystyle{\mathbb{P}(Y_0>0)=1}$. The following holds:
    \begin{itemize}
        \item [(i)] if $\ln(Z_1)$ is integrable, then 
    \begin{equation}
        \mathbb{E}(\ln (Z_1))< 0 \Rightarrow\lim\limits_{n\to \infty}Y_n=0, \hbox{with probability} \ 1 \Rightarrow \mathbb{E}(\ln (Z_1))\leq 0.
    \end{equation}
        \item [(ii)] if $\ln(Z_1)$ is square-integrable, then 
    \begin{equation}
        \lim\limits_{n\to \infty}Y_n=0, \hbox{with probability} \ 1 \Leftrightarrow \mathbb{E}(\ln (Z_1))<0.
    \end{equation}
    \end{itemize}
\end{lemma}

In our case we set $Z_j:=|p_j(z)|$ for a chosen $z\in\mathbb{C}$, where $p_j$ defined as in \eqref{def_p_j}. Recall from \eqref{eq:regs} that 
\begin{equation}
\mathcal{R}^{AS}  = \left\{ z \in \mathbb{C} \colon \lim_{k\to\infty}\left(\left| \eta\right|\cdot \prod_{j=1}^{k} \left| p_j(z) \right|\right) = 0 \text{ with probability } 1 \right\} \nonumber .
\end{equation}
Let us observe that $\left|  p_1(z) \right| = \left| \tau_1\cdot z^2+z+1\right| = \sqrt{f_{a,b}(\tau_1)}$ for $z = a+bi$, $a,b \in\mathbb{R}$, where function $f_{a,b}$ is defined as in Appendix and $\tau_1$ is uniformly distributed over $[0,1]$. From Fact \ref{roots_f_ab} in Appendix it follows that, for all $z\in\mathbb{C}$,  $\mathbb{P}\left( |p_1(z)|>0 \right) = 1$, whereas Fact \ref{sq_int_f} implies that the random variable $\ln\left( \left| p_1(z)\right| \right) = \frac12 \ln\left(  f_{a,b}\left(\tau_1\right) \right)$ is square-integrable. Hence, by Lemma \ref{as_stab_lem_1}(ii) we obtain
\begin{equation} \label{eq:RAS}
\mathcal{R}^{AS}  = \left\{ z \in\mathbb{C}  \colon \mathbb{E}(\ln|p_1(z)|)<0 \right\} .
\end{equation}
\begin{theorem}
\label{prop_AS}
\noindent
    \begin{itemize}
        \item [(i)]  The set $\mathcal{R}^{AS}$ is open and symmetric with respect to the real axis.  
        \item[(ii)] It holds that $\mathcal{I}^{Mid} \subset [-2,0) \subset \mathcal{I}^{AS} \subset \left(-\sqrt{2e},0\right)$.
        \item [(iii)] There exists $r_0\in (0,\infty)$ such that $\mathcal{R}^{AS}\subset \mathbb{C}_{-} \cap \left\{ z\in \mathbb{C} \colon |z|<r_0\} \right\}$.
    \end{itemize}
\end{theorem}
\noindent
{\bf Proof.} In this proof we will refer many times to the family of functions $\left\{ f_{a,b}\colon (a,b)\in\mathbb{R}^2 \right\}$ defined by \eqref{eq:f_formula} and the function $F$ linked to this family via \eqref{def_F}.

Let us notice that $\mathcal{R}^{AS}$ given by \eqref{eq:RAS} is isomorphic with the following set (we will use the same name for both sets):
\begin{equation} \label{RAS_2}
    \mathcal{R}^{AS} = \left\{ (a,b)\in\mathbb{R}^2 \colon F(a,b) < 0 \right\}.
\end{equation}
Since function $F\colon\mathbb{R}^2\to\mathbb{R}$ is continuous (see Proposition \ref{F_cont}), the set $\mathcal{R}^{AS} = F^{-1}\left( (-\infty,0)\right)$ is open. Moreover, $F(a,b) = F(a,-b)$ for all $(a,b)\in\mathbb{R}^2$ (see the proof of Proposition \ref{F_cont}), which implies symmetry of $\mathcal{R}^{AS}$ with respect to the abscissa. This proves (i).

By \eqref{RAS_2},
\begin{equation} \label{IAS}
    \mathcal{I}^{AS} = \left\{ x\in\mathbb{R} \colon F(x,0) < 0 \right\} .
\end{equation}
Fact \ref{fact_F} gives the formula for $F(x,0)$:
\begin{equation}
F(x,0) = \frac{x^2+x+1}{x^2}\ln\left(x^2+x+1\right)-\frac{x+1}{x^2}\ln\left| x+1 \right| - 1    
\end{equation}
for $x\in(-\infty,-1)\cup (-1,0)\cup (0,\infty)$, $F(-1,0)=-1$ and $F(0,0)=0$. 

For $x\in(-1,0)$ we have $1>x^2+x+1>x+1>0$. Function $(0,1) \ni t \mapsto t\ln t - t \in \mathbb{R}$ is decreasing. As a result, 
\begin{equation}
     \left(x^2+x+1\right) \ln\left(x^2+x+1\right) -  \left(x^2+x+1\right) < \left(x+1\right) \ln\left(x+1\right) -  \left(x+1\right),
\end{equation}
which is equivalent to $F(x,0)<0$. We conclude that $(-1,0)\subset \mathcal{I}^{AS}$.

Let us consider a function $g\colon \left( -\infty, -1\right] \to \mathbb{R}$ given by $g(x) = x^2\cdot F(x,0)$, that is
\begin{equation}
    g(x) = \left(x^2+x+1\right) \ln\left(x^2+x+1\right) - \left(x+1\right) \ln\left(-(x+1)\right) - x^2 
\end{equation}
for $x<-1$ and $g(-1)=-1$. Function $g$ is continuous in $(-\infty,-1]$ (since $F$ is continuous) and convex because its second derivative
\begin{equation}
    g''(x) = 2\ln\left(x^2+x+1\right) + \frac{(2x+1)^2}{x^2+x+1}+\frac{1}{-(x+1)}
\end{equation}
is positive for $x\in (-\infty,-1)$ (as each term of the above sum is positive). From Jensen's inequality it follows that
\begin{equation}
    g(x) = g \left( (x+2)\cdot(-1) + (1-x)\cdot (-2)\right) \leq (x+2)\cdot g(-1) + (1-x)\cdot g(-1) < 0
\end{equation}
for all $x\in [-2,-1]$. Hence, $[-2,-1] \subset \mathcal{I}^{AS}$.

We have already shown that $[-2,0) \subset \mathcal{I}^{AS}$. The inclusion $\mathcal{I}^{AS} \subset \left( -\sqrt{2e}, 0 \right)$ follows from the Fact that $\mathcal{R}^{AS} \subset \mathbb{C}^-$, which will be proved later in (iii), and the well-known log sum inequality:
\begin{equation*}
    \sum_{i=1}^n a_i\ln\left(\frac{a_i}{b_i}\right) \geq a\ln\left(\frac{a}{b}\right)
\end{equation*}
for any $n\in\mathbb{Z}_+$ and $a_1,\ldots,a_n,b_1,\ldots,b_n>0$, $a=a_1+\ldots+a_n$, $b=b_1+\ldots,b_n$. In fact, for $x\in \mathcal{I}^{AS} \cap (-\infty,-1)$ we have
\begin{equation}
x^2 > \left(x^2+x+1\right) \ln\left(x^2+x+1\right) - \left(x+1\right) \ln\left(-(x+1)\right) \geq x^2\ln\left(\frac{x^2}{2}\right),
\end{equation}
where the former inequality follows from the condition $F(x,0)<0$ and the latter is log sum inequality with $a_1 = x^2+x+1$, $a_2=-(x+1)$ and $b_1=b_2=1$. We conclude that $1>\ln\left(\frac{x^2}{2}\right)$, which leads to $x>-\sqrt{2e}$ and the proof of (ii) is completed.

Now we will show that $\mathcal{R}^{AS}\subset\mathbb{C}_{-}$. This inclusion is equivalent to the condition $F(a,b) \geq 0$ for all $a\geq 0$ and $b\in\mathbb{R}$. Since the difference
\begin{equation}
\label{diff_f}
    f_{a,b}(t) - f_{0,b}(t) = \left(a^4+2a^2b^2\right)t^2+2\left(a^3+a^2+ab^2\right)t+a^2+2a
\end{equation}
is non-negative for all $t\in [0,1]$, $a\in [0,\infty)$ and $b\in\mathbb{R}$, we have 
\begin{equation}
    F(a,b) = \frac12 \mathbb{E} \left( \ln \left( f_{a,b} (\tau) \right) \right) \geq \frac12 \mathbb{E} \left( \ln \left( f_{0,b} (\tau) \right) \right) = F(0,b).
\end{equation}
Thus, it suffices to show that $F(0,b) \geq 0$ for all $b\in\mathbb{R}$. 

As stated in Fact \ref{fact_F}, $F(0,0)=0$. Hereinafter we assume that $b\neq 0$. Then $\displaystyle f_{0,b} \colon \mathbb{R} \ni t \mapsto b^4t^2-2b^2t+b^2+1 \in \mathbb{R}$ is a quadratic function and its global minimum $b^2$ is achieved for the argument $b^{-2}$. When $|b|\geq 1$, we have $\ln \left( f_{0,b}(t) \right) \geq \ln \left( b^2 \right) \geq \ln 1 = 0$ for all $t\in [0,1]$ and as a result $F(0,b) = \frac12\cdot \mathbb{E} \left( \ln \left( f_{0,b} (\tau) \right) \right) \geq 0$. 

Now we will investigate the remaining case $0<|b|<1$. Recall that $\ln (x) \geq 1-\frac{1}{x}$ for $x>0$,
\begin{equation}
    \text{arctg}(x)-\text{arctg}(y) = \text{arctg}\left(\frac{x-y}{1+xy}\right)
\end{equation} 
for $y<x<0$ and $0<\text{arctg}(x)<x$ for $x>0$. Hence, for $x>0$ it follows that
\begin{equation}
  2F(0,b)  =  \mathbb{E} \left( \ln \left( f_{0,b} (\tau) \right) \right) \geq 1- \mathbb{E} \left( \frac{1}{f_{0,b}(\tau)} \right)
    =  1-\frac{\text{arctg} \left( |b|^3 \right)}{|b|^3} > 0.
\end{equation}
This completes the proof of inclusion $\mathcal{R}^{AS}\subset\mathbb{C}_{-}$.

From (ii) we know that $\mathcal{I}^{AS}$ is bounded. The boundedness of $\mathcal{R}^{AS}$ follows from the following observation: if $z\in\mathcal{R}^{AS}$ and $|z|\geq 4$, then $-|z|\in\mathcal{I}^{AS}$.

Let us consider $z=a+bi\in\mathcal{R}^{AS}$, where $a,b\in\mathbb{R}$, such that $|z|\geq 4$. For $t\in [0,1]$ we have
\begin{align*}
    \left| tz^2+z+1 \right|^2 & - \left| t|z|^2-|z|+1 \right|^2 \\ 
    & = 2 \left( \sqrt{a^2+b^2}+a \right) + 2t \left[ \left(a^2+b^2\right)\sqrt{a^2+b^2}+a\left(a^2+b^2\right)-2b^2 \right] \geq 0.
\end{align*}
We need to provide justification for the last inequality. Firstly, let us observe that $$\sqrt{a^2+b^2}+a \geq  |a|+a \geq 0.$$ Secondly, let us notice that $a<0$ because $z\in\mathcal{R}^{AS}\subset\mathbb{C}_-$ and choose $\alpha\in\mathbb{R}$ such that $b=\alpha \cdot a$. Then $|z|=|a|\sqrt{1+\alpha^2}$ and
    \begin{displaymath}
    \left(a^2+b^2\right)  \sqrt{a^2+b^2}+a\left(a^2+b^2\right)-2b^2  
    = a^2\alpha^2 \left(  |z|\frac{\sqrt{1+\alpha^2}}{\sqrt{1+\alpha^2}+1}-2\right) \geq 0,
\end{displaymath}
since $|z|\geq 4$ and $\frac{\sqrt{1+\alpha^2}}{\sqrt{1+\alpha^2}+1}\geq\frac12$ for all $\alpha\in\mathbb{R}$. 

Hence, for $z\in\mathcal{R}^{AS}$ such that $|z|\geq 4$ the following holds:
\begin{equation*}
z\in\mathcal{R}^{AS} \ \Leftrightarrow \ \mathbb{E} \left( \ln \left| \tau z^2+z+1 \right| \right) <0 \ \Rightarrow \ \mathbb{E} \left( \ln \left| \tau |z|^2-|z|+1 \right| \right) <0 \ \Leftrightarrow \ -|z|\in\mathcal{I}^{AS}.
\end{equation*}
This concludes the proof. \ \ \ $\blacksquare$ 
\subsection{Region $\mathcal{R}^{SP}$ of stability in probability}
Below we prove  analogous result  to Lemma \ref{as_stab_lem_1}, but now we deal with convergence in probability.
\begin{lemma} \label{as_stab_lem_2}
Given a sequence of real-valued, independent and identically distributed random variables $\{Z_n\}_{n\in\mathbb{N}_0}$ with $\displaystyle{\mathbb{P}(Z_1>0)=1}$, consider the sequence of random variables $\{Y_n\}_{n\in\mathbb{N}}$ defined by
    \begin{equation}
        Y_n=\Bigl(\prod\limits_{i=0}^{n-1}Z_i\Bigr)y_0,
    \end{equation}
    where $y_0\in\mathbb{R}_+$. The following holds:
    \begin{itemize}
        \item [(i)] if $\ln(Z_1)$ is integrable, then 
    \begin{equation}
        \mathbb{E}(\ln (Z_1))< 0 \Rightarrow\lim\limits_{n\to \infty}Y_n=0 \ \hbox{in probability} \  \Rightarrow \mathbb{E}(\ln (Z_1))\leq 0.
    \end{equation}
        \item [(ii)] if $\ln(Z_1)$ is square-integrable, then 
    \begin{equation}
        \lim\limits_{n\to \infty}Y_n=0 \ \hbox{in probability} \ \Leftrightarrow \ \mathbb{E}(\ln (Z_1))<0.
    \end{equation}
    \end{itemize}
\end{lemma}
\noindent
{\bf Proof.}
By Lemma \ref{as_stab_lem_1} we obtain 
\begin{equation}
    \mathbb{E}(\ln (Z_1))< 0 \ \ \Rightarrow \ \ Y_n \overset{a.s.}{\longrightarrow} 0 \ \ \Rightarrow \ \ Y_n \overset{\mathbb{P}}{\longrightarrow} 0.
\end{equation}
To prove the second implication in (i), let us suppose that $Y_n \overset{\mathbb{P}}{\longrightarrow} 0$ and $\mathbb{E}(\ln (Z_1)) > 0$. By Riesz theorem, there exists a subsequence $\left( Y_{n_k}\right)_{k=0}^\infty$ of the sequence $\left( Y_n\right)_{n=0}^\infty$ such that $Y_{n_k} \overset{a.s.}{\longrightarrow} 0$. On the other hand, by the strong law of large numbers,
\begin{equation}
    \frac{S_{n_k}}{n_k}=\frac{1}{n_k}\sum_{i=0}^{n_k-1} \ln\left(Z_i\right) \overset{a.s.}{\longrightarrow}\mathbb{E}\left(\ln\left(Z_1\right)\right)>0.
\end{equation}
Thus, $n_k\cdot\frac{S_{n_k}}{n_k} \overset{a.s.}{\longrightarrow} \infty$ and $\displaystyle{
    Y_{n_k} = y_0 \cdot \exp \left( n_k \cdot \frac{S_{n_k}}{n_k} \right)\overset{a.s.}{\longrightarrow} \infty}$. This contradiction proves (i).

To prove part (ii), it suffices to show that the case $Y_n \overset{\mathbb{P}}{\longrightarrow} 0$ and $\mathbb{E}(\ln (Z_1)) = 0$ is impossible. Let us consider this case. Then, by the central limit theorem,
\begin{equation}\label{eq:clt}
    \mathbb{P}\left(  \sum_{i=0}^{n-1} \ln\left(Z_i\right)>\sqrt{n}\cdot\sigma \right) = \mathbb{P}\left( \frac{\frac1{n} \sum\limits_{i=0}^{n-1} \ln\left(Z_i\right)}{\frac{\sigma}{\sqrt{n}}}>1 \right) \xrightarrow{n\to\infty} 1- \Phi(1), 
\end{equation}
where $\sigma =\sqrt{\text{Var}\left(\ln (Z_1)\right)}>0$ and $\Phi$ denotes the CDF of the standard normal distribution. For all $n\in\mathbb{Z}_+$ we have $\sqrt{n}\cdot\sigma>0$. As a result,
\begin{equation}
 \label{eq:convprob}
    \mathbb{P}\left( \sum_{i=0}^{n-1} \ln\left(Z_i\right)>\sqrt{n}\cdot\sigma \right)  \leq \mathbb{P}\left( \sum_{i=0}^{n-1} \ln\left(Z_i\right)> 0 \right) = \mathbb{P}\left( Y_n > y_0 \right) \xrightarrow{n\to\infty} 0,
\end{equation}
because $Y_n \overset{\mathbb{P}}{\longrightarrow} 0$. From \eqref{eq:clt} and \eqref{eq:convprob} it follows that $1-\Phi(1) = 0$, which is a contradiction. Hence, the proof of the lemma is completed. \ \ \ $\blacksquare$\\ \\
From Lemmas \ref{as_stab_lem_1}, \ref{as_stab_lem_2} we get the following.
\begin{corollary} 
\label{cor_as_sp}
Under the assumptions of Lemma \ref{as_stab_lem_2}, if $\ln(Z_1)$ is square-integrable, then
    \begin{equation}
        Y_n \overset{a.s.}{\longrightarrow} 0  \ \Leftrightarrow \ \mathbb{E}(\ln (Z_1))<0 \ \Leftrightarrow \ Y_n \overset{\mathbb{P}}{\longrightarrow} 0.
    \end{equation} 
\end{corollary}
Corollary \ref{cor_as_sp}, \eqref{eq:regs}, and  \eqref{eq:RAS} imply that
\begin{equation} 
\label{RSP}
    \mathcal{R}^{SP}=\{z\in\mathbb{C} \colon \mathbb{E}(\ln|p_1(z)|)<0\}=\mathcal{R}^{AS}.
\end{equation}
Hence, for the randomized Runge-Kutta scheme the notions of asymptotic stability and stability in probability coincide. Furthermore, by \eqref{reg_incl},  \eqref{RSP}, and Theorem \ref{prop_MS} (ii) we have that
\begin{equation}
    \mathcal{R}^{MS}\subset\mathcal{R}^{AS}\cap\mathcal{R}^{Mid}.
\end{equation}
\begin{remark}
The sets $\mathcal{R}^{MS}$, $\mathcal{R}^{AS}$, $\mathcal{R}^{Mid}$ are open (so Borel) and, since they are also bounded, their Lebesgue measure is well defined and finite. Below we present estimates for areas of $\mathcal{R}^{MS}$, $\mathcal{R}^{AS}$ and $\mathcal{R}^{Mid}$
\begin{displaymath}
area(\mathcal{R}^{MS})\approx 3.92\leq area(\mathcal{R}^{AS})=area(\mathcal{R}^{SP})\approx 5.38\leq area(\mathcal{R}^{Mid})\approx 5.87,
\end{displaymath}
however 
\begin{displaymath}
    \mathcal{I}^{MS}\subset\mathcal{I}^{Mid}\subset\mathcal{I}^{AS}=\mathcal{I}^{SP}.
\end{displaymath}
In Figure~\ref{fig:regions} we show the pictures of $\mathcal{R}^{MS}$, $\mathcal{R}^{AS}$, and $\mathcal{R}^{Mid}$ obtained by the Maple package. 
\begin{figure}
\centering
\includegraphics[width=0.7\linewidth]{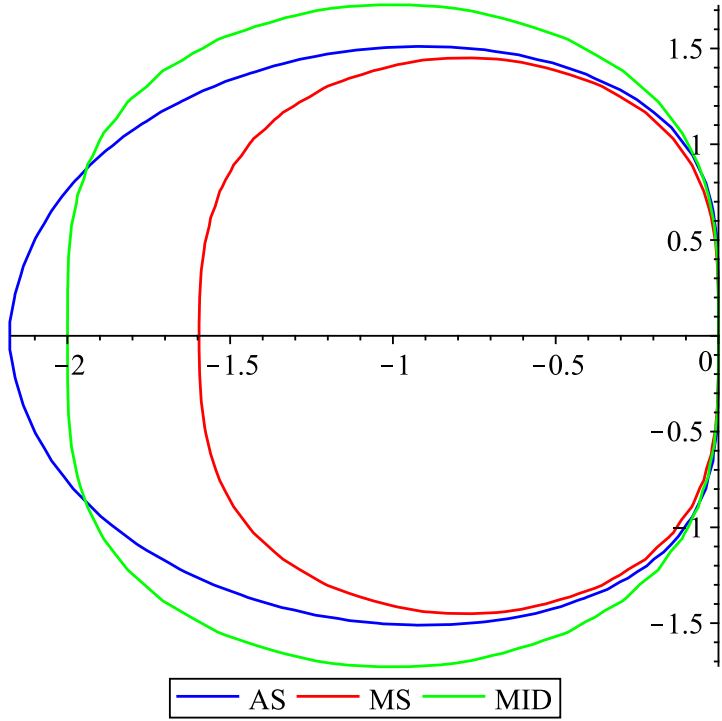}
\caption{Stability regions $\mathcal{R}^{AS}$, $\mathcal{R}^{MS}$ and  $\mathcal{R}^{Mid}$}
\label{fig:regions}\end{figure}
\end{remark}
\begin{remark}
\label{zero_stab_RK}
    Note that \eqref{zero_stab} the randomized Runge-Kutta method under exact information is {\it almost surely $0$-stable} in the sense that there exists $K,h_0>0$ such that for all $f\in F^{\varrho}$, $\varepsilon\in (0,1)$, $\delta\leq\varepsilon$ the following holds
    \begin{equation}
        \max\limits_{0\leq j\leq n}\|V^j-\bar V^j\|\leq K\varepsilon.
    \end{equation}
\end{remark}
\section{Conclusions and future work}
As we have seen, randomization decreases the error (under mild assumptions on right-hand side functions) and, in the case of asymptotic stability, extends the interval of absolute  stability. However, all considered regions of stability are bounded. Therefore, in our future work we intend to consider randomized implicit schemes. We conjecture that at least one of the regions $\mathcal{R}^{MS}$, $\mathcal{R}^{AS}$, $\mathcal{R}^{SP}$ contains $\mathbb{C}_{-}$. \newline\newline
\bf Acknowledgments. \rm This research was partly supported by the National Science Centre, Poland, under project 2017/25/B/ST1/00945.
\section{Appendix}\label{sec:app}

Let us define a function $f_{a,b} \colon \mathbb{R} \to [0,\infty)$ for all pairs $(a,b)\in\mathbb{R}^2$ by the following formula: 
\begin{align}
f_{a,b} (t) & = \left| t(a+bi)^2+a+bi+1 \right|^2 \notag \\ & = \left( a^2+b^2 \right)^2 t^2 + 2 \left( a^2+a^3+ab^2 - b^2 \right) t + (a+1)^2 + b^2. \label{eq:f_formula}
\end{align}
Let us notice that $f_{0,0}\equiv 1$. For each pair $(a,b)\in\mathbb{R}^2\setminus \left\{ (0,0) \right\}$ the function $f_{a,b}$ is quadratic and its discriminant 
\begin{equation}
    \Delta_{a,b} = -4b^2\left( 2a+a^2+b^2\right)^2 \label{eq:delta_ab}
\end{equation}
is non-positive. This leads to Fact \ref{roots_f_ab} below.
\begin{fact} \label{roots_f_ab}
The function $f_{a,b}$ has at most one real root. Moreover,
\ \begin{itemize}
    \item[(i)] there exists $t_0\in\mathbb{R}$ such that $f_{a,b}\left(t_0\right)=0$ if and only if \begin{equation*}
\left\{\begin{matrix}
a \neq 0
\\ b=0
\end{matrix}\right.
\ \ \text{or}  \ \ 
\left\{\begin{matrix}
a \in (-2,0)
\\ b^2 = -a^2-2a
\end{matrix}\right. ,
\end{equation*} 
    \item[(ii)] there exists $t_0\in [0,1]$ such that $f_{a,b}\left(t_0\right)=0$ if and only if \begin{equation*}
\left\{\begin{matrix}
a\in (-\infty,-1]
\\ b=0
\end{matrix}\right. 
\ \ \text{or}  \ \ 
\left\{\begin{matrix}
a \in \left(-2,-\frac12\right]
\\ b^2 = -a^2-2a
\end{matrix}\right. .
\end{equation*}
\end{itemize}
\end{fact}
Let us consider the following function 
\begin{equation} 
\label{def_F}
F \colon \mathbb{R}^2 \ni (a,b) \mapsto \frac12 \mathbb{E}\left( \ln f_{a,b} (\tau)\right) \in \mathbb{R},
\end{equation}
where $\tau$ is a random variable uniformly distributed over the interval $[0,1]$. We will show that $F$ is well-defined and we will express it in explicit form. Let us observe that 
\begin{equation} \label{eq:Fab0}
F (a,b) = \frac12 \int\limits_0^1 \ln \left( A(a,b)t^2+B(a,b)t+C(a,b) \right) \,\mathrm{d}t,
\end{equation}
where $A(a,b) = \left( a^2 + b^2 \right)^2, \ B(a,b) = 2 \left( a^2+a^3+ab^2 - b^2 \right),\  C(a,b)=(a+1)^2+b^2$.
In the case of $a=b=0$ we immediately get $F(0,0)=0$. Hereinafter we assume that at least one of numbers $a,b$ is non-zero, i.e. $a^2+b^2>0$. It implies that $A(a,b)>0$ and $f_{a,b}$ is a quadratic function. We will usually skip arguments $(a,b)$ when using functions $A=A(a,b), B=B(a,b), C=C(a,b)$. The vertex of a parabola $t\mapsto At^2+Bt+C$ has the coordinates $(P,Q)$, where
\begin{align}
P & =P(a,b) = -\frac{B}{2A} = -\frac{a^2+a^3+ab^2-b^2}{\left( a^2+b^2 \right)^2} , \\  Q & =Q(a,b) = \frac{4AC-B^2}{4A} = \frac{b^2\left( 2a+a^2+b^2 \right)^2}{\left( a^2+b^2 \right)^2}.
\end{align} 
The proof of the following fact can be delivered by quite long but straightforward calculations. Hence, we left it to the reader.
\begin{fact}\label{fact_F}
On the axis $b=0$, the function $F$ can be expressed by 
\begin{equation}\label{eq:Fab1}
    F(a,0) = \frac{a^2+a+1}{a^2}\ln\left(a^2+a+1\right)- \frac{a+1}{a^2}\ln\left|a+1\right|-1
\end{equation}
for $a\in \left(-\infty,-1\right) \cup \left(-1,0\right)\cup\left(0,\infty\right)$. Additionally, $F\left( 0,0 \right)=0$ and $F\left( -1,0 \right)=-1$. On the circle with center $(-1,0)$ and radius $1$, function $F$ has the following formula:
\begin{equation}\label{eq:Fab2}
    F(a,b)=\frac{2a+1}{2a}\ln\left|2a+1\right|-1
\end{equation}
for $a\in \left(-2,-\frac12\right)\cup\left(-\frac12,0\right)$ and $b^2=-a^2-2a$. Additionally, $F\left( -\frac12,\frac{\sqrt{3}}{2} \right)=F\left( -\frac12,-\frac{\sqrt{3}}{2} \right)=-1$. For all the remaining pairs $(a,b)\in\mathbb{R}^2$ function F takes the form
\begin{align}\label{eq:Fab3}
F(a,b) & = \frac12 (1-P) \ln (A+B+C) - 1 +\frac{P}{2}\ln C \notag \\ & \ \ \ \ \ \ \ \ \ \ +\sqrt{\frac{Q}{A}}\left[ \text{arctg}\left( (1-P)\sqrt{\frac{A}{Q}} \right) - \text{arctg}\left( -P\sqrt{\frac{A}{Q}} \right) \right].
\end{align}
\end{fact}
\begin{proposition} \label{F_cont}
    The function $F$, defined in \eqref{def_F}, is continuous in $\mathbb{R}^2$.
\end{proposition}
\noindent
{\bf Proof.}
From \eqref{eq:Fab0} it follows that $F(a,b)=F(a,-b)$ for all $(a,b)\in\mathbb{R}^2$, so it suffices to check continuity of function $F$ in $\mathbb{R}\times [0,\infty)$. Let us split $\mathbb{R}\times [0,\infty)$ into the following pairwise disjoint sub-regions associated to equations \eqref{eq:Fab1}--\eqref{eq:Fab3} of function $F$ given in Fact \ref{fact_F}:
\begin{align*}
Z_1 & = \left((-\infty,-1)\cup (-1,0)\cup (0,\infty)\right)\times \{0\}, \\ Z_2 & =  \left\{ (a,b)\in\mathbb{R}\times(0,\infty)\colon a\in \left(-2,-\frac12\right) \cup \left(-\frac12,0\right) \ \wedge \ b^2=-a^2-2a \right\}, \\ Z_3 & =  \left\{ (a,b)\in\mathbb{R}\times (0,\infty) \colon b^2\neq -a^2-2a \right\}.
\end{align*}
Let us notice that $\mathbb{R}\times[0,\infty) = Z_1 \cup Z_2 \cup Z_3 \cup \left\{ X_0, X_1, X_2 \right\}$, where $X_0=(0,0)$, $X_1=(-1,0)$, and $X_2 = \left( -\frac12,\frac{\sqrt{3}}{2} \right)$. Restriction of $F$ to each of the sets $Z_1, \ Z_2, \ Z_3$ is continuous and $Z_3$ is open in $\mathbb{R}\times [0,\infty)$. Hence, it is enough to show that $F$ is continuous in each of the points $X_0, \ X_1, \ X_2$ and in each point belonging to $Z_1$ or $Z_2$. \medskip

Using \eqref{eq:Fab1}, \eqref{eq:Fab2} and the L'H\^{o}pital's rule, it is easy to check that
\begin{equation}
\lim_{(a,b)\to(0,0), \ (a,b)\in Z_1} F(a,b) = \lim_{(a,b)\to(0,0), \ (a,b)\in Z_2} F(a,b) = 0. \label{eq:cont00a}
\end{equation}
Let us notice that $A=\left( a^2 + b^2 \right)^2\to 0$, $B=2 \left( a^2+a^3+ab^2 - b^2 \right)\to 0$ and $C =(a+1)^2+b^2\to 1$, when $(a,b)\to \left( 0,0 \right)$. For $(a,b)\in Z_3$ sufficiently close to $(0,0)$ the following identity holds: 
\begin{equation} \label{eq:ident_arctg}
\text{arctg}\left( (1-P)\sqrt{\frac{A}{Q}} \right) - \text{arctg}\left( -P\sqrt{\frac{A}{Q}} \right) = \text{arctg} \left( \frac{\sqrt{\frac{A}{Q}}}{1-P(1-P)\frac{A}{Q}} \right),
\end{equation}
since $\text{arctg} \,x - \text{arctg}\, y = \text{arctg}\left( \frac{x-y}{1+xy} \right)$ for $x,y\in\mathbb{R}$ such that $\text{arctg}\, x - \text{arctg}\, y \in\left( -\frac{\pi}{2},\frac{\pi}{2} \right)$ and $xy\neq -1$. Let us notice that $\text{arctg}\left( (1-P)\sqrt{\frac{A}{Q}} \right) - \text{arctg}\left( -P\sqrt{\frac{A}{Q}} \right)>0$ because $(1-P)\sqrt{\frac{A}{Q}}>-P\sqrt{\frac{A}{Q}}$. If $P> 1$, then $0>(1-P)\sqrt{\frac{A}{Q}}>-P\sqrt{\frac{A}{Q}}$ and 
\begin{equation*}
\text{arctg}\left( (1-P)\sqrt{\frac{A}{Q}} \right) - \text{arctg}\left( -P\sqrt{\frac{A}{Q}} \right) < \text{arctg}\,0 - \lim\limits_{x\to-\infty}\text{arctg}\,x = \frac{\pi}{2}.
\end{equation*}
For $P < 0$ we proceed analogously. If $P\in [0,1]$, then $B^2 \leq 4A^2$ and as a result
\begin{equation*}
    \frac{A}{Q} = \frac{4A^2}{4AC-B^2} \leq \frac{4A^2}{4AC-4A^2} = \frac{A}{C-A} < 1
\end{equation*}
for $(a,b)\in Z_3$ sufficiently close to $(0,0)$, since $A\to 0$ and $C \to 1$. Thus, $1>(1-P)\sqrt{\frac{A}{Q}}>-P\sqrt{\frac{A}{Q}}>-1$  and
\begin{equation*}
\text{arctg}\left( (1-P)\sqrt{\frac{A}{Q}} \right) - \text{arctg}\left( -P\sqrt{\frac{A}{Q}} \right) < \text{arctg}\,1 - \text{arctg}(-1) =\frac{\pi}{2}.
\end{equation*}
for $(a,b)\in Z_3$ sufficiently close to $(0,0)$. Given that $ \frac{1}{1-P(1-P)\frac{A}{Q}} =\frac{2C}{2C+B} - \frac{B^2}{4A}\cdot \frac{2}{2C+B}$, 
we obtain
\begin{align}
\sqrt{\frac{Q}{A}}\cdot\text{arctg} \left( \frac{\sqrt{\frac{A}{Q}}}{1-P(1-P)\frac{A}{Q}} \right) & = \frac{1}{1-P(1-P)\frac{A}{Q}} \cdot \frac{\text{arctg} \left( \frac{\sqrt{\frac{A}{Q}}}{1-P(1-P)\frac{A}{Q}} \right)}{ \frac{\sqrt{\frac{A}{Q}}}{1-P(1-P)\frac{A}{Q}} } \notag \\ & = \left(C - \frac{B^2}{4A}\right)\cdot g_1(a,b), \label{eq:cont00_1}
\end{align}
where 
\begin{equation}
g_1(a,b)=\frac{2}{2C+B}\cdot \frac{\text{arctg} \left( \frac{\sqrt{\frac{A}{Q}}}{1-P(1-P)\frac{A}{Q}} \right)}{ \frac{\sqrt{\frac{A}{Q}}}{1-P(1-P)\frac{A}{Q}} } \to 1 \ \ \text{as} \ (a,b) \in Z_3 \ \text{and} \ (a,b)\to(0,0)
\end{equation}
because $\frac{\sqrt{\frac{A}{Q}}}{1-P(1-P)\frac{A}{Q}} = \frac{\sqrt{4AC-B^2}}{2C+B} \to 0$. Furthermore,
\begin{align}
\frac12 (1-P)  \ln (A+B+C) + \frac{P}{2}\ln C & = \frac12 \ln (A+B+C) +\frac{B}{4A}\cdot\frac{A+B}{C} \ln\left( 1 + \frac{A+B}{C}\right)^{\frac{C}{A+B}} \notag \\ & = g_2(a,b) +\frac{B^2}{4A} \cdot g_3(a,b),  \label{eq:cont00_2}
\end{align}
where 
\begin{equation}
g_2(a,b) =  \frac12 \ln (A+B+C)  + \frac{B}{4C} \cdot \ln\left(1 + \frac{A+B}{C} \right)^{\frac{C}{A+B}} \to 0
\end{equation}
and
\begin{equation} \label{eq:cont00_g3}
g_3(a,b) = \frac{1}{C} \cdot \ln\left( 1 + \frac{A+B}{C}\right)^{\frac{C}{A+B}} \to 1
\end{equation}
as $(a,b)\in Z_3$ and $(a,b)\to \left( 0,0 \right)$. We can show that $\frac{B^2}{4A}$ is bounded for $0<a^2+b^2<1$. For this purpose let us introduce polar coordinates: $a=r\cos\varphi$ and $b=r\sin\varphi$, where $r\in (0,1), \ \varphi\in[0,2\pi)$. Then
\begin{equation} \label{eq:cont00_bound}
\frac{B^2}{4A} = \left( r\cos\varphi+\cos 2\varphi \right)^2 \leq \left( r\left|\cos\varphi\right|+\left|\cos 2\varphi\right| \right)^2 \leq (r+1)^2 < 4.
\end{equation}
By \eqref{eq:Fab3} and \eqref{eq:ident_arctg}--\eqref{eq:cont00_bound}:
\begin{equation*}
F(a,b) = g_2(a,b)+\frac{B^2}{4A} \cdot \left( g_3(a,b)-g_1(a,b) \right) + C\cdot g_1(a,b)-1 \to 0
\end{equation*}
as $(a,b)\in Z_3$ and $(a,b)\to \left( 0,0 \right)$. This combined with \eqref{eq:cont00a} implies continuity of $F$ in $(0,0)$. \medskip

Now we check the continuity of function $F$ in point $X_1=(-1,0)$. From \eqref{eq:Fab1} it is easy to see that $F(a,b) \to -1$ as $(a,b) \in Z_1$ and $(a,b)\to X_1$. Let us notice that $A \to 1$, $B \to 0$, $C \to 0$, $P \to 0$, $Q\to 0$ and $\frac12 (1-P) \ln (A+B+C)  \to 0$, when $(a,b)\in Z_3$ tends to $X_1$. Moreover, $\frac{P}{2}\ln C \to 0$ because 
\begin{displaymath}
\lim_{(a,b)\to(-1,0)} \frac{B\ln C}{2}  = 
\lim_{(a,b)\to(0,0)} \left[ (a-2)\left(a^2+b^2\right)\ln\left( a^2+b^2 \right)+a \ln\left( a^2+b^2 \right) \right] = 0.
\end{displaymath}
Furthermore,
\begin{align*}
& \sqrt{\frac{Q}{A}} \cdot  \left[ \text{arctg}\left( (1-P)\sqrt{\frac{A}{Q}} \right) - \text{arctg}\left( -P\sqrt{\frac{A}{Q}} \right) \right] \to 0,
\end{align*}
when $(a,b)\in Z_3$ and $(a,b)\to X_1$, since $\frac{Q}{A} \to 0$ and arctangent is a bounded function. As a result, $F(a,b)\to -1$, when $(a,b)\in Z_3$ tends to $X_1$ and continuity of $F$ in $X_1$ follows. \medskip

We check the continuity of function $F$ in point $X_2=\left( -\frac12,\frac{\sqrt{3}}{2} \right)$. Let us observe that $A\to 1$, $B\to -2$, $C\to 1$, $P=-\frac{B}{2A} \to 1$ and $Q=\frac{4AC-B^2}{4A} \to 0$ when $(a,b)\in Z_3$ and $(a,b)\to X_2$. Thus,
\begin{align}
&\frac{P}{2}\ln C \to 0 \ \text{ and } \ \sqrt{\frac{Q}{A}} \cdot  \left[ \text{arctg}\left( (1-P)\sqrt{\frac{A}{Q}} \right) - \text{arctg}\left( -P\sqrt{\frac{A}{Q}} \right) \right] \to 0,\label{eq:X2_1}
\end{align}
when $(a,b)\in Z_3$ and $(a,b)\to X_2$. Moreover, since $\sqrt{A+B+C}\cdot \ln \sqrt{A+B+C} \to 0$ and $\frac{A+\frac{B}{2}} {\sqrt{A+B+C}}$ is bounded for $(a,b)\in Z_3$ in some neighbourhood of $X_2$, we get
\begin{equation} \label{eq:X2_2}
\frac12 (1-P) \ln (A+B+C) = \frac{1}{A} \cdot \frac{A+\frac{B}{2}} {\sqrt{A+B+C}}\cdot\sqrt{A+B+C}\cdot \ln \sqrt{A+B+C} \to 0
\end{equation}
as $(a,b)\in Z_3$ tends to $X_2$. From \eqref{eq:Fab3}, \eqref{eq:X2_1} and \eqref{eq:X2_2} it follows that $F(a,b)\to -1 = F\left(X_2 \right)$, when $(a,b)\to X_2$ and $(a,b)\in Z_3$. Thus, $F$ is continuous in $X_2$. Since the boundedness of $\frac{A+\frac{B}{2}}{\sqrt{A+B+C}}$ is not straightforward, we will provide a justification. To analyse this expression, it will be convenient to use the polar coordinates:
\begin{equation*}
\left\{ \begin{matrix} a=r\cos\varphi-\frac12 \\ b=r\sin\varphi+\frac{\sqrt{3}}{2} \end{matrix} \right. \ \text{ with } \ r> 0, \ \varphi\in [0,2\pi).
\end{equation*}
Then
\begin{equation*}
\frac{A+\frac{B}{2}}{\sqrt{A+B+C}} = \resizebox{0.8\hsize}{!}{ $\frac{ \textstyle{r^3 + r^2\left( 2\sqrt{3}\sin\varphi-\cos\varphi \right)+r\left( \frac52+\sin^2\varphi-\sqrt{3}\sin\varphi\cos\varphi \right)}+\frac{\sqrt{3}}{2}\sin\varphi-\frac32\cos\varphi }{\sqrt{r^2+r\cdot 2\sqrt{3}\sin\varphi+3}} $ }
\end{equation*}
and we can observe that the above expression is bounded for $r\leq\frac12$:
\begin{align*}
\frac{\left|A+\frac{B}{2}\right|}{\sqrt{A+B+C}} & \leq \resizebox{0.8\hsize}{!}{ $\frac{ \textstyle{r^3 + r^2\left( 2\sqrt{3}\left|\sin\varphi\right|+\left|\cos\varphi\right| \right)+r\left( \frac52+\sin^2\varphi+\sqrt{3}\left|\sin\varphi\cos\varphi\right| \right)}+\frac{\sqrt{3}}{2}\left|\sin\varphi\right|+\frac32\left|\cos\varphi\right| }{\sqrt{r^2-r\cdot 2\sqrt{3}\left|\sin\varphi\right|+3}}$} \notag \\ & \leq \frac{ \frac18 + \frac14\left( 2\sqrt{3}+1 \right)+\frac12\left( \frac52+1+\sqrt{3} \right)+\frac{\sqrt{3}}{2}+\frac32}{\sqrt{3-\sqrt{3}}}.
\end{align*} \smallskip

We check the continuity of function $F$ in points from $Z_1$. Let us consider $\left( a_0, b_0 \right) \in Z_1$. Then $a_0 \in (-\infty,-1) \cup (-1,0) \cup (0,\infty)$ and $b_0=0$. Let us observe that $A \to a_0^4$, $B\to 2a^2_0(a_0+1)$, $C \to (a_0+1)^2$, $P \to -\frac{a_0+1}{a_0^2}$ and $Q \to 0$, when $(a,b)\in Z_3$ and $(a,b)\to (a_0,0)$. Hence,
\begin{equation*}
\frac12 (1-P) \ln (A+B+C)  \to  \frac{a_0^2+a_0+1}{a_0^2} \ln \left(a_0^2+a_0+1\right)
\end{equation*}
and 
\begin{equation*}
\frac{P}{2}\ln C \to  -\frac{a_0+1}{a_0^2} \ln \left|a_0+1 \right|.
\end{equation*}
Since $\frac{Q}{A}\to 0$,  we obtain
\begin{equation*}
\sqrt{\frac{Q}{A}} \cdot  \left[ \text{arctg}\left( (1-P)\sqrt{\frac{A}{Q}} \right) - \text{arctg}\left( -P\sqrt{\frac{A}{Q}} \right) \right]\to 0.
\end{equation*}
We combine the above considerations with \eqref{eq:Fab3}, which results in the following:
\begin{equation*}
F(a,b) \to \frac{a_0^2+a_0+1}{a_0^2}\ln\left(a_0^2+a_0+1\right)- \frac{a_0+1}{a_0^2}\ln|a_0+1|-1 = F\left(a_0,0\right),
\end{equation*}
when $(a,b)\in Z_3$ and $(a,b)\to \left( a_0,0 \right)$. This means that $F$ is continuous in $\left( a_0,0\right)$ if $a_0 \neq -2$. The point $(-2,0)$ is special because in each its punctured neighbourhood there are points not only from $Z_1$ and $Z_3$, but also from $Z_2$. Hence, we need to calculate the following limit:
\begin{equation*}
\lim_{(a,b)\to(-2,0), \ (a,b)\in Z_2} F(a,b)  = \lim_{a\to -2^+ } \left[  \frac{2a+1}{2a}\ln\left(-2a-1\right)-1 \right] = \frac34 \ln 3 -1 = F(-2,0)
\end{equation*}
and only now we can conclude that $F$ is continuous in $(-2,0)$. \medskip

We check the continuity of function $F$ in points from $Z_2$. Let us consider $\left( a_0, b_0 \right) \in Z_2$. Then $a_0 \in \left(-2,-\frac12\right) \cup \left(-\frac12,0\right)$ and $b^2_0=-a_0^2-2a_0$. Let us observe that $A\to 4a_0^2$, $B\to 4a_0$, $C\to 1$, $P\to -\frac{1}{2a_0}$ and $Q\to 0$, when $(a,b)\in Z_3$ and $(a,b)\to (a_0,b_0)$. As a result, by \eqref{eq:Fab3} we obtain
\begin{equation*}
F(a,b) \to \frac{2a_0+1}{2a_0}\ln\left|2a_0+1\right|-1 = F\left(a_0,b_0\right),
\end{equation*}
when $(a,b)\in Z_3$ and $(a,b)\to \left( a_0,b_0 \right)$. Therefore, $F$ is continuous in $\left( a_0,b_0\right)$. 

Finally, we can conclude that the function $F$ is continuous in $\mathbb{R}^2$. \ \ \ $\blacksquare$
\begin{fact}
\label{sq_int_f}
The random variable $\ln f_{a,b} (\tau)$, where $\tau$ is uniformly distributed over the interval $[0,1]$, is square integrable for all $(a,b)\in\mathbb{R}^2$.
\end{fact}
\noindent
{\bf Proof.} 
As stated in Fact \ref{roots_f_ab}, the function $[0,1]\ni t \mapsto \left(\ln f_{a,b}(t)\right)^2 \in\mathbb{R}$ is continuous for all pairs $(a,b)\in\mathbb{R}^2$ but those satisfying one of the following conditions:
\begin{equation*}
    1^\circ \ \ \left\{\begin{matrix}
a\in (-\infty,-1]
\\ b=0
\end{matrix}\right. 
\ \ \ \text{or}  \ \ \ 2^\circ \ \ 
\left\{\begin{matrix}
a \in \left(-2,-\frac12\right]
\\ b^2 = -a^2-2a
\end{matrix}\right. .
\end{equation*}
In case $1^\circ$ we obtain
\begin{displaymath}
\mathbb{E}\left( \ln f_{a,b} (\tau) \right)^2=\int\limits_0^1 \left(\ln f_{a,b}(t)\right)^2 \,\mathrm{d}t = \frac{4}{a^2}\left[ G(-a-1)+G\left(a^2+a+1\right) - 2G(0) \right] < \infty,
\end{displaymath}
where $G(x)=x(\ln x)^2-2x\ln x+2x$ is continuous on $(0,\infty)$ and can be continuously extended on $[0,\infty)$ with $G(0)=0$. In case $2^\circ$ we proceed similarly and arrive at
$$\int\limits_0^1 \left(\ln f_{a,b}(t)\right)^2 \,\mathrm{d}t  = -\frac{1}{2a}\left[ G(1)+G\left(-2a-1\right) - 2G(0) \right] < \infty.$$ 
This completes the proof. \ \ \ $\blacksquare$

\end{document}